\documentclass[a4paper,10pt,reqno]{amsart}

\usepackage{amsmath,graphics}
\usepackage{amssymb}
\usepackage{amsfonts}
\usepackage{latexsym}
\usepackage[mathscr]{eucal}
\usepackage[dvips]{graphicx}

\theoremstyle{plain}
\newtheorem{thm}{Theorem}[section]
\newtheorem{propo}[thm]{Proposition}
\newtheorem{lem}[thm]{Lemma}
\newtheorem{cor}[thm]{Corollary}
\newtheorem{defi}[thm]{Definition}

\theoremstyle{definition}

\newtheorem{rem}[thm]{Remark}

\renewcommand{\Re}{{\rm Re}}
\renewcommand{\Im}{{\rm Im}}
\newcommand{\R}{\mathbb{R}}
\newcommand{\C}{\mathbb{C}}
\newcommand{\Z}{\mathbb{Z}}
\newcommand{\N}{\mathbb{N}}

\renewcommand{\H}{{\mathbb{H}^2}}

\newcommand{\rr}{\varrho}
\newcommand{\lt}{\mathcal{L}}
\newcommand{\eb}{\mathbf{e}}
\newcommand{\GL}{\rm{GL}}
\newcommand{\Id}{\rm{Id}}
\newcommand{\tr}{\rm{tr}}
\newcommand{\rank}{\rm{rank}}

\title[]{On the spectrum of twisted Laplacians and the Teichm\"uller representation}

\author[F.~Naud]{Fr\'ed\'eric Naud}
\address{%
Institut Math\'ematique de Jussieu\\
Universit\'e Pierre et Marie Curie, 4 place Jussieu, 75252 Paris Cedex 05\\
France.
}
\email{frederic.naud@imj-prg.fr}
\author[P.~Spilioti]{Polyxeni Spilioti}
\address{ Department of Mathematics\\
Ny Munkegade 118, building 1530, 4238000 Aarhus C\\
Denmark
}
\email{spilioti@math.au.dk}

\subjclass{}

\keywords{}

\begin{document}
\bibliographystyle{plain}

\maketitle

\begin{abstract}
Given a compact hyperbolic surface $X=\Gamma \backslash \H$ and a linear non-unitary representation $\rr:\Gamma\rightarrow \mathrm{GL}(V)$, we investigate the spectrum of
the twisted Laplacian $\Delta_\rr$ acting on sections of the associated flat vector bundle $E_\rr\rightarrow X$. We show that this non self-adjoint operator has its spectrum inside
a parabola related to a critical exponent $\delta$ of the representation. In the case where $\rr$ is of Teichm\"uller type
we then exhibit explicit parabolic regions determined by another constant $\delta_0$ with an asymptotic spectral density which improves Weyl's law by a power factor. Both $\delta$ and $\delta_0$
are uniquely determined by the so-called ``Manhattan curve'' related to the representation $\rr$.
\end{abstract}

\bigskip \tableofcontents

\section{Introduction and main results}
The Laplace-Beltrami operator is a key object in Riemannian geometry: its spectrum encodes some fine properties of the geometry of the manifold and the dynamics of the geodesic flow. This operator is self-adjoint and has discrete spectrum on compact manifolds. It has natural generalizations to vector bundles and remains self-adjoint if the vector bundle comes equipped with a natural hermitian metric. In particular, if this vector bundle is flat, the holonomy of a loop depends only on its homotopy class and there is an associated natural unitary representation of the fundamental group. Conversely, any unitary finite dimensional representation of the fundamental group of a compact manifold give rise to a flat vector bundle on which a natural 
self-adjoint Bochner-Laplace operator lives. There are however natural situations where non-unitary representations of the fundamental group lead to consider an associated non-self-adjoint Laplacian, whose discrete spectrum remains rather mysterious. In this paper, we will study non-self-adjoint Bochner-Laplace operators associated to linear, non-unitary representations of surface groups. There is a very rich litterature on this subject which leads to a wide class of non-self-adjoint Laplacians for which, apart from the general results pertaining to elliptic operators, almost nothing is known.

More precisely, let $\H$ denote the hyperbolic plane and let $X=\Gamma \backslash \H$ be a compact hyperbolic surface obtained by quotienting $\H$ by a discrete co-compact subgroup of $\mathrm{PSL}_2(\R)$.
Given a linear representation $\rr:\Gamma\rightarrow \mathrm{GL}(V)$, where $V$ is a complex $d$-dimensional vector space, we can easily build a flat vector Bundle $E_\rr$ over $X$ by considering the quotient
of $\H\times V$ by the $\Gamma$-action defined by $\gamma.(z, v)\colon=(\gamma(z), \rr(\gamma) v)$. Smooth sections of $E_\rr$ are then identified with elements $F\in C^\infty(\H,V)$ satisfying
$$\forall\ \gamma \in \Gamma,\ F(\gamma z)=\rr(\gamma)F(z).$$
Let $\Delta_\rr$ be the usual hyperbolic Laplacian acting diagonally on smooth $\rr$-equivariant sections, then this operator can be identified with the twisted Bochner-Laplace operator 
$$\Delta_\rr=-\mathrm{tr}( (\nabla^{E_\rr})^2),$$
where $\nabla^{E_\rr}$ is a flat connection on $E_\rr$. Choosing a hermitian metric on $V$, one can define $L^2(X,E_\rr)$ as the space of square integrable sections of $E_\rr$. The operator $\Delta_\rr$ is in general not self-adjoint.
However, it has a self-adjoint
principal symbol and hence it has nice spectral properties (\cite{shubin1987pseudodifferential}).
It is an elliptic operator and, because of the compactness of the base $X$, it has a unique closed extension with dense domain in $L^2(X,E_\rr)$, with a compact resolvent. Its spectrum is therefore discrete
and moreover is contained in a translate of a positive cone $C\subset \C$ with $\R^{+}\subset\C     $ (\cite{shubin1987pseudodifferential},   \cite[Lemma 2.1]{Muller1}).
In addition, the spectrum of the twisted Laplacian satisfies a Weyl law, as showed by M\"uller in a pioneering work \cite{Muller1}. He also proved a general trace formula for integral operators induced  by the twisted Laplacian,
using Paley-Wiener test functions.
This work generalizes Selberg's theory in a non-self-adjoint setting for all locally symmetric space of split rank one. The higher rank case
is considered by Shen, where he used the M\"uller's Selberg trace formula to study the Ruelle zeta function and its relation to the complex analytic torsion (\cite{Shen}). In general, the twisted Laplacians are defined as coupling operators to a flat vector bundle
(\cite[Section 2]{cappell2010complex}). 
One can define then also the twisted Hodge Laplacian, acting on the space of vector-valued differential forms. Consequently, 
there exists an extension of the usual notion of the analytic torsion (\cite{ray1971r}) to the case of non self-adjoint Hodge Laplacians. This is the Cappell-Miller torsion (\cite{cappell2010complex}) and the refined analytic torsion of Braverman and Kappeler (\cite{braverman2007refined},\cite{braverman2008refined}).

\bigskip
Let now $\mathcal{P}$ denote the set of primitive conjugacy classes in $\Gamma$ and if $\gamma \in \mathcal{P}$, let $\ell(\gamma)$ be the length of the associated primitive closed geodesic on $X$. We define the critical exponent $\delta$ of the representation $\rr$ by 
$$\delta(\rr):= \inf \left \{ s>0\ :\ \sum_{\gamma\in \mathcal{P},k\geq 1} \vert \mathrm{tr}(\rr(\gamma^k))\vert e^{-sk\ell(\gamma)}<\infty \right \}.$$
A simple argument using the fact that the word length and the displacement function (with respect to hyperbolic distance in $\H$) are equivalent shows that there exist
$M,C>0$ and a set of representatives $\{ \gamma \}$ of conjugacy classes such that, $\vert \mathrm{tr}(\rr(\gamma))\vert\leq Me^{C\ell(\gamma)}$, hence $\delta(\rr)$ is finite. Since $\rr$ is assumed to be non unitary, we will have in general $\delta(\rr)>1$. For certain type of representations $\rr$, this critical exponent is related to the celebrated Manhattan curve, see below for details.

We then define the parabolic regions $\mathcal{C}_\sigma\subset \C$ for all $\sigma \in (1/2,+\infty)$ by
$$\mathcal{C}_\sigma:=\left \{z=x+iy \in \C\ :\ x\geq \sigma(1-\sigma)+\frac{y^2}{(1-2\sigma)^2}\right \}.$$ 
We will also use the notation $P_\sigma:=\partial \mathcal{C}_\sigma$ to denote the parabola given by the equation $\{x+iy\ :\ x=\sigma(1-\sigma)+\frac{y^2}{(1-2\sigma)^2}\}$. In the degenerate case where $\sigma=1/2$, we set
$P_\sigma=\mathcal{C}_\sigma=[1/4,+\infty)$.
Our first result is the following.
\begin{thm}
\label{main1} Let $\Delta_\rr$ be the twisted Laplacian defined as above, then one always have 
$$\mathrm{Sp}(\Delta_\rr)\subset \mathcal{C}_{\delta(\rr)}.$$
\end{thm}
We will show in $\S 3$ that this inclusion is in some sense sharp: there exists families of linear representations $\rr$ such that $\Delta_\rr$ satisfies 
$$\inf_{\lambda \in \mathrm{Sp}(\Delta_\rr)} \Re(\lambda) =\delta(1-\delta).$$
It is tempting to believe that one always has this identity, but this is not true: in $\S 3$ we exhibit a family of linear representations for which $\delta(1-\delta) \not \in \mathrm{Sp}(\Delta_\rr)$.

The proof of Theorem \ref{main1} is standard and follows from the already known properties of the twisted Selberg
zeta functions and a Heat trace argument. As proved by Muller in \cite{Muller1}, the twisted Laplacian satisfies a Weyl law: as $T\rightarrow +\infty$,
$$\#\{ \lambda \in  \mathrm{Sp}(\Delta_\rr)\ :\ \vert \lambda \vert \leq T\}=\frac{\mathrm{dim}(V)\mathrm{Vol}(X)}{4\pi}T+o(T).$$

We prove a much more precise fact for some specific type of representations which we call ``of Teichm\"uller type". Given
a fixed hyperbolic, compact surface $X=\Gamma\backslash \H$, we will view the {\it Teichm\"uller space} $\mathcal{T}(X)$ as the space  of {\it faithfull discrete representations }
$$\rho:\Gamma \rightarrow \mathrm{PSL}_2(\R),$$
quotiented by conjugations in $\mathrm{PSL}_2(\R)$.
We denote the hyperbolic distance on $\H$ by $d(z,w)$, we will use the upper half-plane model $\H=\{\Im(z)>0\}$ in the definition below.
\begin{defi} 
A linear representation $\rr:\Gamma\rightarrow \mathrm{GL}(V)$ is said to be of {\it Teichm\"uller type} if there exist a norm $\Vert . \Vert_0$ on $\mathrm{GL}(V)$, a representation
$\rho \in \mathcal{T}(X)$ and constants $C,\beta>0$ such that for all $\gamma \in \Gamma$,
\begin{enumerate}
\item $\Vert \rr(\gamma) \Vert_0 \leq C e^{\beta d(i,\rho(\gamma)i)},$
\item $ C^{-1} e^{\beta \ell(\rho(\gamma))}\leq \vert \mathrm{tr}(\rr(\gamma))\vert\leq C e^{\beta \ell(\rho(\gamma))}.$
\end{enumerate}
\end{defi}

 Let $\rr:\Gamma \rightarrow \mathrm{GL}(V_\rr)$
  be a representation of Teichm\"uller type, as defined above. Clearly we have
 $$\delta(\rr)=\inf \left \{s>0\ :\ \sum_{\gamma,k} e^{\beta k\ell_\rho(\gamma)-sk\ell(\gamma)}<\infty \right \}, $$
 where $\rho:\Gamma \rightarrow \mathrm{PSL}_2(\R)$ is a Teichm\"uller deformation of $\Gamma$ and $\ell_\rho(\gamma)$ stands for $\ell(\rho(\gamma))$. Following the work of Burger \cite{Burger},
 let $\mathcal{M}_\rho$ be the convex region defined by 
 $$\mathcal{M}_\rho:=\left \{ (a,b)\in \R^2\ :\  \sum_{\gamma,k} e^{-ak\ell_\rho(\gamma)-bk\ell(\gamma)}<\infty \right \}.$$
 Then, the boundary of $\mathcal{M}_\rho$ is a convex graph $\mathfrak{C}$ which contains the points $(0,1)$ and $(1,0)$, called ''the Manhattan curve". This curve is a straight line if and only if $\Gamma$ and $\rho(\Gamma)$ are conjugated in $\mathrm{PSL}_2(\R)$, see \cite{Burger,Sharp}.
 It was shown by Sharp in \cite{Sharp} that this curve is actually real-analytic. In our setting, $\delta(\rr)$ is the unique real number such that $(-\beta,\delta(\rr))$ belongs to $\mathfrak{C}$, see the figure below.
\begin{center}
\includegraphics[scale=0.3]{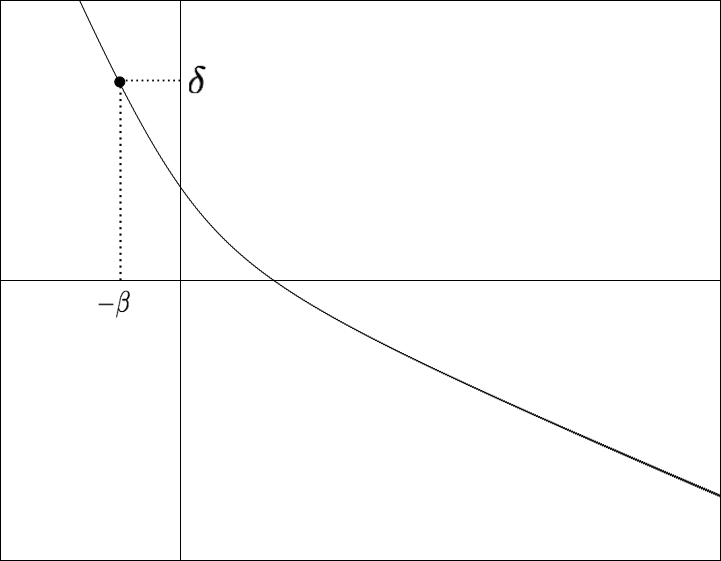}\\
{\tiny The Manhattan curve $\mathfrak{C}$ and the exponent $\delta(\rr)$. }
\end{center}

The main result of this paper is the following improved Weyl law.
\begin{thm}
 \label{main2}
 Assume that $\rr$ is a a linear representation of Teichm\"uller type. There exists $\delta/2<\delta_0<\delta$ such that for all $\delta_0<\sigma\leq \delta$,  one can find
 $\epsilon(\sigma)>0$ such that as $r\rightarrow +\infty$, we have 
 $$M_\sigma(r):=\#\{ \lambda \in \mathrm{Sp}(\Delta_\rr)\setminus \mathcal{C}_\sigma\ :\ r\leq \vert \lambda \vert \leq r +\sqrt{r}\}=O(r^{1/2-\epsilon(\sigma)}).$$
 Moreover, similarly as $\delta$, the constant $\delta_0$ is uniquely determined by $\beta$ and the Manhattan curve: if $(-2\beta,\eta)\in \mathfrak{C}$, then $\delta_0=\eta/2$.
\end{thm}
We have the following corollary.
\begin{cor} For all $\sigma>\delta_0$, 
as $T\rightarrow +\infty$ we have
$$N_\sigma(T):=\#\{ \lambda \in \mathrm{Sp}(\Delta_\rr)\setminus \mathcal{C}_\sigma\ :\ \vert \lambda \vert \leq T\}=O(T^{1-\epsilon(\sigma)}),$$
which beats the Weyl law by a polynomial factor. As a consequence, Weyl law actually holds in $\mathcal{C}_{\sigma}$ for all $\sigma>\delta_0$: as $T\rightarrow +\infty$,
$$\#\{ \lambda \in  \mathrm{Sp}(\Delta_\rr)\cap \mathcal{C}_{\sigma}\ :\ \vert \lambda \vert \leq T\}=\frac{\mathrm{dim}(V)\mathrm{Vol}(X)}{4\pi}T+o(T).$$
\end{cor}
\noindent {\it Proof}. Consider the sequence $(n_p)_{p\geq 0}$ defined by $n_0:=1$ and $n_{p+1}=n_p+\sqrt{n_p}$ for all $p\geq 1$. It is a standard fact that as $p\rightarrow +\infty$ we have 
$$n_p=\frac{p^2}{4}(1+o(1)).$$
We fix $\sigma>\delta_0$, let $T>>1$ and choose the smallest integer $p$ such that $T\leq n_{p+1}$. We have $p=O(\sqrt{T})$, while $n_p<T$. We then write
$$N_\sigma(T)\leq O(1)+\sum_{k=0}^p M_\sigma(n_k)\leq  O(1)+(p+1)\max_{0 \leq k \leq p}M_\sigma (n_k)=O\left(T^{\frac{1}{2}+\frac{1}{2}-\epsilon(\sigma)}\right),$$
and the proof is done. $\square$

\bigskip
The above theorem therefore shows that in the ``high frequency'' regime when we have $\Re(\lambda)\rightarrow +\infty$, the bulk of the spectrum lies in a smaller parabolic region $\mathcal{C}_{\delta_0}$. 
\begin{center}
\includegraphics[scale=0.3]{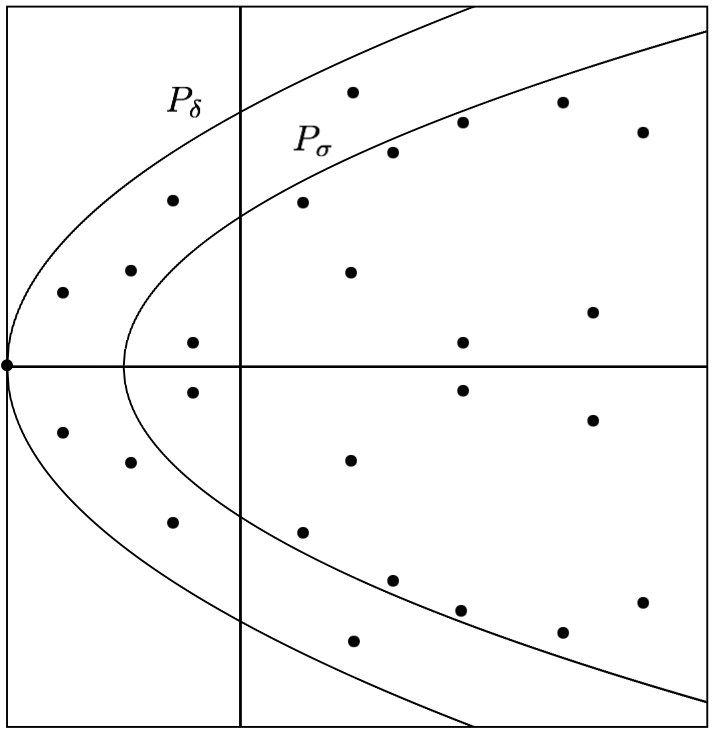}\\
{\tiny Spectrum of $\Delta_\rr$ with $P_\delta$ and $P_{\sigma}$, the bulk of the spectrum being inside $\mathcal{C}_{\sigma}$.}
\end{center}
Some explicit bounds on $\delta_0$ involving the ``geodesic strech" can be derived, see $\S 6$. The computable examples of $\S 3$ show that $\delta_0$ is sharp.

In the case when $\rr$ is a {\it one-dimensional non-unitary representation},  Anantharaman \cite{Anantharaman1} proves a similar result as Theorem \ref{main2}. Indeed, if we have
$$\rr(\gamma)=e^{\int_{\widetilde{\gamma}} \omega},$$
where $\omega$ is a harmonic real-valued $1$-form on $X=\Gamma \backslash \H$, and $\widetilde{\gamma}\in H^1(X,\Z)$ is the homology class determined by $\gamma$. The spectrum of $\Delta_\rr$ acting on sections of the line bundle $E_\rr$ is the same as the spectrum of
$$\Delta_\omega f:=\Delta_X f+2\omega(\nabla f),$$
acting on functions on $X$. This spectral problem can therefore be understood in terms of a generalized stationary damped wave equation on $X$, and the semi-classical technique of Sj\"ostrand \cite{Sjostrand1} can be refined as done in \cite{Anantharaman1} to obtain the one-dimensional analog of Theorem \ref{main2}. Unfortunately, this one-dimensional technique cannot be applied for higher dimensional representations. Hence our method is fairly different, but our result can be viewed as a generalization of \cite{Anantharaman1} to the higher dimensional representations case. It is also interesting  to observe that in \cite{Anantharaman1}, independently from \cite{Muller1}, a trace formula for $\Delta_\rho$ is derived and utilized, by an analytic continuation argument of the ``usual'' Selberg trace formula for a unitary representation.

\bigskip
The paper is organized as follows. In $\S 2$, we recall the known facts on twisted Laplacians, including the heat trace formula and the main result on the Selberg zeta function from \cite{SF1}. We then prove easily Theorem \ref{main1}. In $\S 3$, we exhibit representations of Teichm\"uller type and compute the spectrum for several examples: depending on the parity of the dimension of the representation space, these spectra lie on parabolas and are related to (self-adjoint) Dirac and Laplace spectra. In $\S 4$, we recall the so-called Bowen-series coding for co-compact Fuchsian groups and relate $\delta$ to the topological pressure by using some basic Teichm\"uller theory and thermodynamical formalism. The proof of Theorem \ref{main2} really starts in $\S 6$, where we relate the counting problem to an estimate from above of certain Hilbert-Schmidt norms for a family of  transfer operators. These transfer operators act on a family of (vector valued) Bergman spaces of holomorphic functions and 
rather delicate off-diagonal estimates for Bergman kernels of shrinking domains are needed: they are proved in the Appendix. We finally analyse $\delta_0$ and prove some explicit bounds via Thermodynamical formalism. At the technical level, the core of the proof of Theorem \ref{main2} bears some similarities with \cite{JN1,Naud1}, but there are some serious additional difficulties and differences which make it interesting: the non-unitarity of the representation involved, and the need to use more general complex domains (hence the non-trivial Bounds on Bergman kernels).

\bigskip
{\bf Acknowledgements.} We thank Julien March\'e and Andres Sambarino for several discussions around the examples of $\S 3$ and pointing out
several possible generalizations. We also thank Maxime Wolff for pointing out references on Thurston's asymetric distance.

\section{Twisted Selberg zeta functions and the main parabola}
In this section we prove Theorem \ref{main1}, which is completely general. 
We define first the twisted Bochner--Laplace operator as introduced  in \cite[Section 4]{Muller1}. This operator acts on sections of a vector bundle over $X$. It is an elliptic operator which is in general not self-adjoint. Nevertheless, it has a self-adjoint principle symbol (with respect to the choice of a metric) and hence it has qualitatively similar spectral properties.

Let $\rr$ be a finite-dimensional, complex representation
\begin{equation*}
\rr\colon\Gamma\rightarrow \GL(V_{\rr})\end{equation*}
of $\Gamma$.
Denote by $E_{\rr}=V_\rr\times_\Gamma\H\rightarrow X$ the associated flat vector bundle over $X$, equipped with a flat connection $\nabla^{E_{\rr}}$.
We recall the construction of the twisted Bochner--Laplace operator $\Delta_\rr^\sharp$, acting on smooth sections of $E_{\rr}$.
The second covariant derivative $(\nabla^{E_\rr})^2$ is defined by
$$
(\nabla^{E_\rr})^2_{V,W}:=\nabla_{V}^{E_\rr}\nabla_{W}^{E_\rr}-\nabla^{E_\rr}_{\nabla^{LC}_{V}W}, 
$$
where $V,W\in C^{\infty}(X,TX)$, $TX$ is the tangent bundle of $X$,
and $\nabla^{LC}$ denotes the Levi--Civita connection on $TX$. 
The twisted Bochner--Laplace operator $\Delta_\rr^\sharp$ is defined to be the corresponding connection Laplacian on $E_\rr$, i.e., the negative of the trace of the second covariant derivative:
\begin{equation}\label{eq:DefDeltaSharp}
\Delta_\rr^{\sharp}:= -\tr\big((\nabla^{E_\rr})^2\big).
\end{equation}

Locally, this operator is described as follows.
We consider an open subset $U$ of $X$ such that $E_{\rr}\lvert_{U}$ is trivial, 
i.e.,  $E_\rr\lvert_U\cong U\times\C^m$, where $m=\rank(E_\rr)=\dim V_{\rr}$.
Let $e_1,\ldots,e_m$ be any basis of flat sections of $E_{\rr}\lvert_{U}$. Then, each 
$\phi \in C^{\infty}(U,E_{\rr}\lvert_{U})$ can be written as
$
\phi=\sum_{i=1}^{m}\phi_{i} \otimes e_{i},
$
where $\phi_{i}\in C^{\infty}(U)$, $i=1,\ldots,m$. Then, 
\begin{equation}
\nabla_{Y}^{E_{\rr}}\phi=\sum_{i=1}^{m}\nabla_{Y}\phi_{i}\otimes e_{i}, \qquad (Y\in C^{\infty}(X,TX)).\label{eq:ConnectionLocally}
\end{equation}
The local expression above is independent of the choice of the basis of flat sections of $E_{\rr}\lvert_{U}$, since the transition maps comparing flat sections are constant. By {\eqref{eq:DefDeltaSharp}} and \eqref{eq:ConnectionLocally}, the twisted Bochner--Laplace operator acting on $C^{\infty}(U,E_{\chi}\lvert_{U})$ is given by
\begin{equation}\label{sharploc}
\Delta_\rr^\sharp\phi=\sum_{i=1}^{m}(\Delta\phi_{i})\otimes e_{i},
\end{equation}
where $\Delta$ denotes the usual Laplace--Beltrami operator on $X$ for the hyperbolic metric.
Let now $\widetilde{E}_{\rr}$ be the pullback to $\widetilde{X}=\H$ of $E_{\rr}$.
Then,
\begin{equation*}
\widetilde{E}_{\rr}\cong \widetilde{X}\times V_{\rr},
\end{equation*}
and
\begin{equation}\label{iso}
C^{\infty}(\widetilde{X},\widetilde{E}_{\rr})\cong  C^{\infty}(\widetilde{X})\otimes V_{\rr}.
\end{equation}
With respect to the isomorphism {\eqref{iso}}, it follows from {\eqref{sharploc}} that the lift $\widetilde{\Delta}^{\sharp}_\rr$ of $\Delta_\rr^{\sharp}$ to $\widetilde{E}_{\rr}$ takes the form
\begin{equation}\label{eq:DeltaSharpTilde}
\widetilde{\Delta}^{\sharp}_\rr = \widetilde{\Delta}\otimes\Id_{V_{\rr}},
\end{equation}
where $\widetilde{\Delta}$ is the Laplace--Beltrami operator on $\H$.

If we choose a Hermitian metric on $E_{\rr}$, then $\Delta_\rr^{\sharp}$ acts in $L^{2}(X,E_{\rr})$ with domain $C^{\infty}(X,E_\rr)$. However, it is not a formally self-adjoint operator in general. By {\eqref{sharploc}}, $\Delta_\rr^{\sharp}$ has principal symbol 
\begin{equation*}
\sigma_{\Delta_{\rr}^{\sharp}}(x,\xi)=\lVert \xi \rVert^ {2}_{x} \Id_{({E_{\rr})_{x}}}
\qquad (x\in X, \xi\in T_{x}^{*}X).
\end{equation*}
Hence, $\Delta_\rr^\sharp$ is an elliptic, second order differential 
operator with the following spectral properties:
its spectrum is discrete and contained in a translate of a positive cone $C\subset \C$ such that $\R^{+}\subset C$.
This fact follows from classical spectral theory of elliptic operators, under the assumption of the compactness of the manifold.
We refer the reader to \cite{shubin1987pseudodifferential}, and also \cite[Lemma 2.1]{Muller1}. From now on, for notational simplicity, we will denote this Laplacian by $\Delta_\rr$.

The semi-group $e^{-t\Delta_\rr}$ is well defined for $t>0$ and trace class on $L^2(X,E_\rr)$. The following ``heat trace formula" \cite{SF1} holds:
\begin{equation}
 \label{heattrace}
 \mathrm{tr}(e^{-t\Delta_\rr})=\sum_{\mu \in  \mathrm{Sp}(\Delta_\rr)} e^{-t\mu}
\end{equation}
$$=\frac{1}{4\pi^2} \mathrm{dim}(V_\rr)\mathrm{Vol}(X)\int_\R e^{-t(\lambda^2+1/4)}\lambda \pi \tanh(\lambda \pi)d\lambda$$
$$ +\frac{1}{2\pi}\sum_{\gamma \in \mathcal{P}, k\geq 1} \mathrm{tr}(\rr(\gamma^k))\frac{\ell(\gamma)}{2\sinh(k\ell(\gamma)/2)} e^{-t/4}e^{-\frac{(k\ell(\gamma))^2}{4t}}.$$
By working with this heat trace one can then deduce a ``resolvent trace formula"  \cite{SF1} and obtain the following fact on the twisted Selberg zeta function.
\begin{thm}
\label{zetadiv}
For $\Re(s)$ large enough, set
$$Z_\Gamma(s,\rr):=\prod_{k\geq 0} \prod_{\gamma \in \mathcal{P}} \det \left( I-\rr(\gamma)e^{-(s+k)\ell(\gamma)}\right).$$
Then $Z_\Gamma(s,\rr)$ has a holomorphic extension to $\C$. Its zeros set is the divisor (with multiplicities)
$$\bigcup_j \{1/2\pm i r_j \}\cup \bigcup_{k=0}^\infty (2g-2)\mathrm{dim}(V_\rr)(2k+1)\{-k\},$$
where $\{ 1/4+r_j^2\}$ is the spectrum of $\Delta_\rr$ and eigenvalues are repeated according to their multiplicity. In addition the following functional equation holds for all $s\in \C\setminus \Z$,
$$Z_\Gamma(s,\rr)=Z_\Gamma(1-s,\rr)\mathrm{exp}\left(\mathrm{dim}(V_\rr)\mathrm{Vol}(X)\int_0^{s-1/2}  r\tan(\pi r)dr \right).$$
\end{thm}

\begin{rem}
 By Theorem 2.1, the twisted Selberg zeta function $Z_\Gamma(s,\rr)$ has a holomorphic extension to $\C$. Indeed, by considering the ``resolvent trace formula" \cite[eq. (4.21)]{SF1}, we have that the logarithmic derivative $L(s,\rr)$ of $Z_\Gamma(s,\rr)$ has simple poles with positive residues. 
This follows from Proposition 4.2.4 and Proposition 4.2.5 in \cite{SF1}. Still, we can not conclude more information about the order of the zeros of $Z_\Gamma(s,\rr)$, since the ``spectral'' and the ``trivial'' zeros could overlap.
\end{rem}

\subsection{Proof of Theorem \ref{main1}}
For $\Re(s)$ large, we have by a direct calculation
$$\frac{Z'_\Gamma(s,\rr)}{Z_\Gamma(s,\rr)}=\sum_{\gamma \in \mathcal{P}, k\geq 1} \frac{\ell(\gamma)}{1-e^{-k\ell(\gamma)}} \mathrm{tr}(\rr(\gamma^k))e^{-sk\ell(\gamma)}.$$
By definition of $\delta(\rr)$, we know that the right hand side converges absolutely for $\Re(s)>\delta$, which shows that $Z_\Gamma(s,\rr)$ has no zeros in the half-plane $\{\Re(s)>\delta \}$.
Using the functional equation, we deduce that non-trivial zeros of $Z_\Gamma(s,\rr)$ are in the strip $\{1-\delta\leq \Re(s)\leq \delta\}$. Given $s_j(1-s_j)\in \mathrm{Sp}(\Delta_\rr)$ with
$s_j=1/2\pm ir_j$, then we know that $s_j$ must be a zero of $Z_\Gamma(s,\rr)$: therefore we get that
$$\mathrm{Sp}(\Delta_\rr)\subset \Z^-\cup \mathcal{C}_\delta.$$
We have to discard the possible negative integer eigenvalues that are below $\delta(1-\delta)$. To this end we go back to the heat trace formula (\ref{heattrace}) in the regime $t\rightarrow +\infty$.
Let us denote by $N_\rr(x)$ the counting function
$$N_\rr(x):=\sum_{k\ell(\gamma)\leq x} \vert \mathrm{tr}(\rr(\gamma^k))\vert. $$
By definition of $\delta$, we have the crude bound valid for all $\epsilon>0$ as $x\rightarrow \infty$:
$$N_\rr(x)=O_\epsilon( e^{(\delta+\epsilon)x} ).$$
Set $\alpha:=\inf_{\mu \in \mathrm{Sp}(\Delta_\rr)} \Re(\mu)$. By discreteness of the spectrum, we know that there are finitely many eigenvalues $\mu$ such that $\Re(\mu)=\delta$, and therefore a simple argument involving Lebesgue dominated
convergence shows that there exists $C_0>0$ such that for all $t$ large,
$$C_0 e^{-t\alpha}\leq \vert \sum_j e^{-t\mu_j} \vert=\vert \mathrm{tr}(e^{-t\Delta_\rr} )\vert. $$
Using the heat trace formula we get for $t$ large
$$C_0 e^{-t\alpha}\leq O(e^{-t/4})+O\left ( e^{-t/4}\sum_{\gamma \in \mathcal{P}, k\geq 1} \vert \mathrm{tr}(\rr(\gamma^k))\vert \frac{\ell(\gamma)}{2\sinh(k\ell(\gamma)/2)} e^{-\frac{(k\ell(\gamma))^2}{4t}}    \right).$$
Furthermore, we have 
$$ S(t):=\sum_{\gamma \in \mathcal{P}, k\geq 1} \vert \mathrm{tr}(\rr(\gamma^k))\vert \frac{\ell(\gamma)}{2\sinh(k\ell(\gamma)/2)} e^{-\frac{(k\ell(\gamma))^2}{4t}}$$
$$ =O\left ( \int_0^\infty xe^{-x/2}e^{-\frac{x^2}{4t}}dN_\rr (x)\right).$$
A Stieltjes integration by parts combined with the crude bound on $N_\rr(x)$ shows that for all $\epsilon>0$ we have :
$$S(t)=O\left (  \int_0^\infty e^{(\delta+\epsilon-1/2)x}e^{-\frac{x^2}{4t}}dx\right).$$
 It is a standard calculation (Laplace transform of a Gaussian) to show that as $a\rightarrow \infty$ we have
 $$\int_0^\infty e^{a v}e^{-v^2}dv=e^{a^2/4}(\sqrt{\pi}+o(1)).$$
 Using this asymptotic and a change of variable show that for all $\epsilon>0$ and all $t$ large,
 $$ S(t)=O\left ( \sqrt{t} e^{(-\delta(1-\delta)+1/4+\epsilon)t}   \right ),$$
 which implies that as $t\rightarrow + \infty$, 
 $$C_0 e^{-\alpha t}\leq O\left ( \sqrt{t} e^{(-\delta(1-\delta)+\epsilon)t}   \right )$$
 and therefore $\alpha\geq \delta(1-\delta)$. The proof is done. $\square$
 
 \subsection{On the bottom of the spectrum}
From Theorem \ref{main1}, we know that $$\inf_{\mu \in \mathrm{Sp}(\Delta_\rr)} \Re(\mu)\geq \delta(1-\delta).$$
 In $\S 3$, we will give examples for which this inequality is strict. There is however a large class of representations for which
the equality holds true. We prove the following.
\begin{propo}
Assume that for all $\gamma \in \Gamma$ we have $\mathrm{tr}(\rr(\gamma))\geq 0$. Then we have
$$\inf_{\mu \in \mathrm{Sp}(\Delta_\rr)} \Re(\mu)=\delta(1-\delta).$$
\end{propo}
\noindent {\it Proof}. Assume that the inequality is strict, then we know that the meromorphic function on $\C$
$$s\mapsto\frac{Z_\Gamma'(s,\rr)}{Z_\Gamma(s,\rr)}$$
has no pole at $s=\delta$, i.e. is analytic at $s=\delta$. On the other hand we have for all $\Re(s)>\delta$
$$\frac{Z'_\Gamma(s,\rr)}{Z_\Gamma(s,\rr)}=\sum_{\gamma \in \mathcal{P}, k\geq 1} \frac{\ell(\gamma)}{1-e^{-k\ell(\gamma)}} \mathrm{tr}(\rr(\gamma^k))e^{-sk\ell(\gamma)},$$
which is a Dirichlet series of the type $\sum_n a_n e^{-s\lambda_n}$ with $a_n\geq 0$, whose abcissa of absolute convergence is precisely $\delta$. A classical theorem of Landau \cite{Landau} tells us that $s=\delta$
must be an analytic singularity, a contradiction. $\square$

We remark that if $\rr$ is a representation of Teichm\"uller type, then we must have $\delta(\rr)>1$, as depicted on the Manhattan curve. Therefore $\delta(1-\delta)<0$ which means that $\Delta_\rr$ has indeed negative eigenvalues
under the hypothesis of positivity of the character $\mathrm{tr}(\rr(\gamma))$. 

If the character is {\it real-valued}, then from the heat trace formula one deduces that for all $t>0$ we have
$$\sum_j e^{-t\mu_j}=\sum_j e^{-t\overline{\mu_j}},$$
which implies that $\mathrm{Sp}(\Delta_\rr)=\overline{\mathrm{Sp}(\Delta_\rr)}$, i.e., the spectrum is symmetric with respect to the $x$-axis. Alternatively, we can recover this fact from the identity
$$Z_\Gamma(s,\rr)=\overline{Z_\Gamma(\overline{s},\rr)},$$
which follows from the equality (valid for $\Re(s)$ large)
$$Z_\Gamma(s,\rr)=\exp \left (-\sum_\gamma \sum_{n\geq 1} \frac{\mathrm{tr}(\rr(\gamma^n))}{n}\frac{e^{-s\ell(\gamma^n)}}{1-e^{-\ell(\gamma^n)}}  \right),$$
and uniqueness of analytic continuation.

\section{Examples of linear representations of surface groups} 
In this section, we give examples of representations $\rr : \Gamma \rightarrow \mathrm{GL}(V)$ which are of {\it Teichm\"uller type}. For some of them the spectrum of $\Delta_\rr$ can be explicitely computed in term of Laplace or Dirac spectrum on the base surface $X$. We then make some remarks about the bottom of the spectrum of $\Delta_\rr$ for these specific representations.
\subsection{The adjoint representation}
Let $G=\mathrm{SL}_2(\R)$ and denote by $V=\mathfrak g$ the Lie algebra of $G$, seen as a real vector space spanned by $Y_0,Y_1,Y_2$, where 
$$Y_0=\left ( \begin{array}{cc}0&1\\0&0  \end{array}\right),\ Y_1=\left ( \begin{array}{cc}0&0\\1&0  \end{array}\right),\ Y_2=\left ( \begin{array}{cc}1&0\\0&-1  \end{array}\right).$$
The adjoint representation $\mathrm{Ad}:G\rightarrow \mathrm{GL}(V)$ is defined by $\mathrm{Ad}(g)(X)=gXg^{-1}$, and leaves invariant the killing form on $V$ given by
$$B(X,Y)=\mathrm{tr}(XY).$$
Notice that $\mathrm{Ad}(-I)=I$, hence $\mathrm{Ad}$ induces a representation $$r_a:\mathrm{PSL}_2(\R)\rightarrow \mathrm{GL}(V).$$
In the basis $(Y_0,Y_1,Y_2)$, the associated quadratic form to $B$ is 
$$B(U,U)=2u_2^2+2u_0u_1$$
for $U=u_0Y_0+u_1Y_1+u_2 Y_2$, which is a quadratic form of signature $(2,1)$. This means that in a suitable basis $(X_0,X_1,X_2)$, this quadratic form can be expressed as
$$B(X,X)=-x_0^2+x_1^2+x_2^2,$$
for $X=x_0X_0+x_1X_1+x_2X_2$.
The Hyperboloid model of the hyperbolic plane $\H$ is then
$$\H:=\{ X\in V\ :\ B(X,X)=-1\ \mathrm{and}\ x_0>0\}.$$
It is a standard fact that $\mathrm{Ad}$ induces an isomorphism between $PSL_2(\R)$ and $SO^+(2,1)$ which is the group of direct isometries of the hyperboloid model. Notice that hyperbolic distance in the hyperboloid model
is given by $$\cosh(d(U,V))=-B(U,V).$$ In particular, if $\mathrm{Ad}(g)$ is given in the basis $(X_0,X_1,X_2)$ by the matrix
$$M=\left (\begin{array}{ccc} a&b&c\\ d&e&f\\ g&h&i     \end{array}  \right), $$
then the fact that $M$ and $M^{-1}$ preserve the killing form imply that we have the formula
$$\Vert M \Vert^2=4a^2-1=4\cosh^2d(\mathrm{Ad}(g)o,o)-1, $$
where $o=(1,0,0)\in \H$. Given an hyperbolic element $\gamma \in \mathrm{PSL}_2(\R)$, we can always conjugate it to 
$$\widetilde{\gamma} =\left ( \begin{array}{cc} e^{\ell(\gamma)/2}&0\\ 0& e^{-\ell(\gamma/2}\end{array}   \right),$$
so that a simple computation in the basis $(Y_0,Y_1,Y_2)$ yields
\begin{equation}
\label{vapadj}
\mathrm{Ad(\widetilde{\gamma})}Y_0=e^{\ell(\gamma)}Y_0,\  \mathrm{Ad(\widetilde{\gamma})}Y_1=e^{-\ell(\gamma)}Y_1,\ \mathrm{Ad(\widetilde{\gamma})}Y_2=Y_2,
\end{equation}
\begin{equation}
\label{traceadj}
\mathrm{tr}(\mathrm{Ad}(\gamma))=\mathrm{tr}(\mathrm{Ad}(\widetilde{\gamma}))= 1+e^{\ell(\gamma)}+e^{-\ell(\gamma)}.
\end{equation}
All the above computations show that for all $\rho \in \mathcal{T}(\Gamma)$, $r_a\circ \rho:\Gamma \rightarrow SO^+(2,1)$ is indeed a Teichm\"uller type representation with $\beta=1$. 

In the specific case $\rr=r_a\circ id$, then we have
by the prime orbit theorem $\delta(\rr)=2$. The spectrum of $\Delta_\rr$ is actually computable. Indeed, using formulas (\ref{vapadj}), we can express the twisted Selberg zeta function as follows. For all $\Re(s)$ large,
we have
$$Z_\Gamma(s,\rr)=\prod_{k\in \N}\prod_{\gamma \in \mathcal{P}} \det \left( I-\rr(\gamma)e^{-(s+k)\ell(\gamma)}\right)$$
$$=\prod_{k\in \N}\prod_{\gamma \in \mathcal{P}}\left(1-e^{-(s+k-1)\ell(\gamma)}\right)\left(1-e^{-(s+k+1)\ell(\gamma)}\right)\left(1-e^{-(s+k)\ell(\gamma)}\right)$$
$$=Z_\Gamma(s-1)Z_\Gamma(s+1)Z_\Gamma(s),$$
where $Z_\Gamma(s)$ is the non-twisted Selberg zeta function of $X=\Gamma \backslash \H$. The divisor of this zeta function is (with multiplicities)
$$\{1/2\pm i t_j\}\cup \bigcup_{k=0}^\infty (2g-2)(2k+1)\{-k\},$$
where the spectrum of the Laplacian $\Delta$ on the base surface $X$ is given by
$$\mathrm{Sp}(\Delta)=\{ 1/4+t_j^2,\ j\in \N \},$$
see for example Hejhal \cite{Hejhal}.
We point out  that $1$ is a simple zero of $Z_\Gamma(s)$.
The divisor of $Z_\Gamma(s,\rr)$ is therefore
$$\left \{\frac{1}{2}\pm i t_j \right \}\cup\left \{\frac{3}{2}\pm i t_j \right\}\cup \left \{-\frac{1}{2}\pm i t_j \right \}\cup 4(2g-2)\{0\}$$
$$\cup (2g-2)\{1\}\cup \bigcup_{k=1}^\infty (2g-2)3(2k+1)\{-k\}.$$
By comparing it with the known divisor of $Z_\Gamma(s,\rr)$ from Theorem \ref{zetadiv}, we deduce that
the bottom of the spectrum of $\Delta_\rr$ is therefore $-2=\delta(1-\delta)$, and is a simple eigenvalue. All but finitely many eigenvalues of $\Delta_\rr$ are 
included in $\mathcal{C}_{3/2}$. Notice that we have in this specific example
$\mathrm{Sp}(\Delta_\rr)\subset \R \cup P_{3/2}$.
\begin{center}
\includegraphics[scale=0.3]{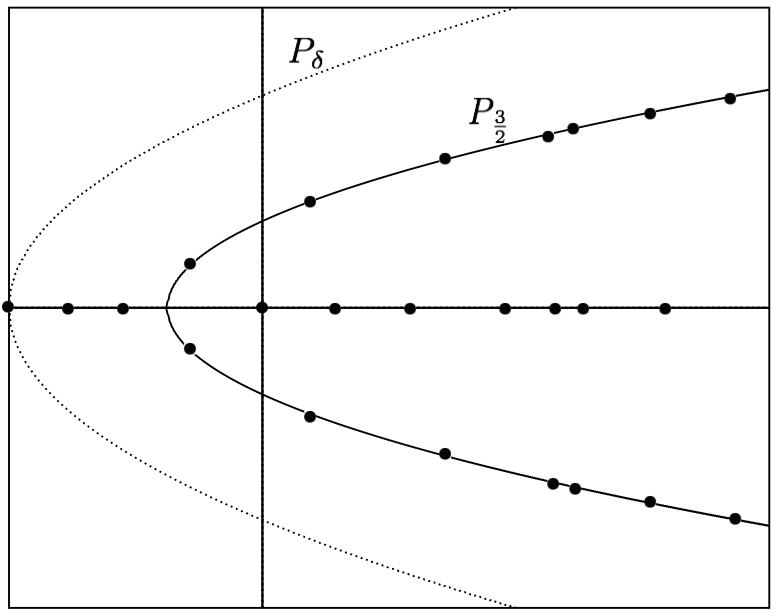}\\
{\tiny Spectrum of $\Delta_\rr$ with $P_\delta=\partial \mathcal{C}_\delta$ and $P_{3/2}=\partial \mathcal{C}_{3/2}$. }
\end{center}

\subsection{Spin structures and $\mathrm{SL}_2(\R)$-lifts}
Let us denote by $$\pi:\mathrm{SL}_2(\R)\rightarrow \mathrm{PSL}_2(\R)$$ the natural projection and set $\widetilde{\Gamma}=\pi^{-1}(\Gamma)$, which is a subgroup of $\mathrm{SL}_2(\R)$ with $-I\in \widetilde{\Gamma}$.
Given a character $\chi:\widetilde{\Gamma}\rightarrow \{+1,-1\}$, with $\chi(-I)=-1$, we can define a linear representation $r_\chi:\Gamma \rightarrow \mathrm{SL}_2(\R)$ by setting
$$r_\chi(\gamma):=\chi(\widetilde{\gamma})\widetilde{\gamma},$$
where $\widetilde{\gamma}$ is such that $\pi(\widetilde{\gamma})=\gamma$. Clearly $r_\chi(\gamma)$ does not depend on the choice of $\widetilde{\gamma}$ and it is straightforward to check that this is a homomorphism.
The existence of such a $\chi$ is non trivial and corresponds to the existence of spin structures on the riemann surface $X$, there are $2^{2g}$ such choices, where $g$ is the genus of $X$, see \cite{Atiyah}. To such a spin structure one can associate
a complex line bundle $T^{1/2}$ on $X$ which is a square root of the canonical line bundle $T$ (bundle of holomorphic $1$-forms on $X$), i.e. $T^{1/2}\otimes T^{1/2}=T$.
Spinors (sections of $T^{1/2}$) correspond, from the point of view of automorphic forms, to smooth functions $f:\H\rightarrow \C$ such that for all $\gamma \in \Gamma$,
$$f(\gamma z)=\chi(\widetilde{\gamma}) \left ( \frac{cz+d}{\vert cz+d\vert}   \right)f(z), $$
where $\widetilde{\gamma}=\left ( \begin{array}{cc}a&b\\c&d   \end{array}    \right)$. 
The space of such automorphic forms is denoted by $S(1)$. The Laplacian $\Delta_1$ is given in the Poincar\'e model by 
$$\Delta_1=-y^2 \left(\frac{\partial^2}{\partial x^2}+ \frac{\partial^2}{\partial y^2}\right)+iy\frac{\partial}{\partial x},$$
and leaves invariant $S(1)$, see Hejhal \cite{Hejhal}, chapter 4. This operator $\Delta_1$, is essentially self-adjoint, has a positive discrete spectrum starting at $1/4$, see \cite{Hejhal} and references herein.
Following Roelcke, Hejhal has developped a Selberg trace formula for $S(1)$, and there is a natural associated Selberg zeta function $Z_\Gamma(s,\chi)$ which is defined as follows for $\Re(s)>1$:
$$Z_\Gamma(s,\chi)=\prod_{\gamma \in \mathcal{P}}\prod_{k\geq 0}\left (1-\chi(\widetilde{\gamma})e^{-(s+k)\ell(\gamma)}     \right),$$
where $\widetilde{\gamma}\in \mathrm{SL}_2(\R)$ is such that $\pi(\widetilde{\gamma})=\gamma$ and $\mathrm{tr}(\widetilde{\gamma})>0$. It is important to stress that $\chi$ is not a character of $\Gamma$ and that the sign changes given by
$\chi(\widetilde{\gamma})$ depend on the homology class of a lift of the closed geodesic $\gamma$ to the unit tangent bundle $T^1X$. By a trace formula argument, see \cite{Sarnak1}, 
we know that $Z_\Gamma(s,\chi)$ has a holomorphic extension to $\C$ and its trivial zeros are located at $s=-\frac{1}{2}-n$, $n\in \N$ with multiplicities $(2n+2)(2g-2)$. The non-trivial zeros are given by $s_j=\frac{1}{2}\pm i r_j$, where 
the spectrum of $\Delta_1$ is $\{1/4+r_j^2,\ j\in \N\}$. Since 
the bottom of the spectrum of $\Delta_1$ starts at $1/4$, all non-trivial zeros of $Z_\Gamma(s,\chi)$ actually lie on the line $\{\Re(s)=1/2\}$.

There is also an interpretation of $r_j$ as spectrum of a genuine Dirac operator acting on a rank $2$ vector bundle, see \cite{Bolte1}.
As in the previous example, it is definitely possible to express the spectrum of $\Delta_{r_\chi}$ in terms of the spectrum of $\Delta_1$. Indeed, we have for all $\gamma \in \Gamma$,
$$r_\chi(\gamma)=\chi(\widetilde{\gamma})\widetilde{\gamma},$$
where $\widetilde{\gamma}$ can be chosen with $\mathrm{tr}(\widetilde{\gamma})>0$ so that $\widetilde{\gamma}$ is conjugate in $\mathrm{SL}_2(\R)$ to 
$$\left ( \begin{array}{cc}e^{\ell(\gamma)/2}&0\\0& e^{-\ell(\gamma)/2}  \end{array}    \right).$$
A similar computation as before shows that for $\Re(s)$ large we have
$$Z_\Gamma(s,r_\chi)=Z_\Gamma(s+1/2,\chi)Z_\Gamma(s-1/2,\chi).$$

The divisor of $Z_\Gamma(s,r_\chi)$ is therefore 
$$ \{ 1\pm ir_j\}\cup \{\pm i r_j \}\cup \bigcup_{k=0}^{\infty} (2g-2)2(2k+1)\{-k\}.$$
 By comparing with the known divisor from Theorem \ref{zetadiv}, we therefore conclude that the spectrum of $\Delta_{r_\chi}$ is included in the parabola $P_0$, see picture below.
\begin{center}
\includegraphics[scale=0.3]{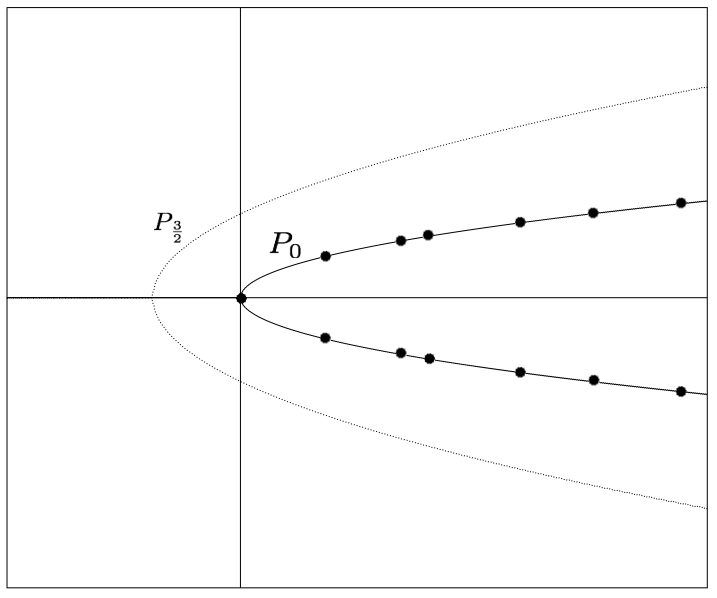}\\
{\tiny Spectrum of $\Delta_{r_\chi}$ and Dirac spectrum on $P_0$. }
\end{center}
By the prime orbit theorem, we know that $\delta(r_\chi)=3/2$, and we point out that $\delta(1-\delta)=-3/4$ 
while the bottom of the spectrum of $\Delta_{r_\chi}$ is $0$ and has multiplicity $2\mathrm{dim}( h^0(X,E^{1/2}))$, where $h^0(X,E^{1/2})$ denotes the space of holomorphic sections of $E^{1/2}$, see \cite{Hejhal}, chapter 4.
Given a Teichm\"uller representation $\rho:\Gamma \rightarrow \mathrm{PSL}_2(\R)$, we can choose a spin structure for $\rho(\Gamma)$ and obtain a more general linear representation $$\rr=r_\chi \circ \rho:\Gamma \rightarrow \mathrm{SL}_2(\R).$$
It is straightforward to check that this is a representation of Teichm\"uller type with $\beta=1/2$, since for all 
$$\gamma \simeq \left ( \begin{array}{cc}a&b\\c&d   \end{array}    \right)\in \mathrm{PSL}_2(\R),$$ the following formula holds in the Poincar\'e half-plane model:
$$a^2+b^2+c^2+d^2=2\cosh(d(i,\gamma i)). $$
\subsection{The Fuchsian representation in $\mathrm{SL}_{n+1}(\R)$}
Let $\mathcal{H}_n$ denote the space of bivariate homogeneous polynomials of degree $n$, that is 
$$\mathcal{H}_n=\mathrm{Span}\{ P_n^j(x,y):=x^jy^{n-j},\ j=0,\ldots,n\}. $$
Any element $M=\left ( \begin{array}{cc}a&b\\c&d   \end{array}    \right) \in \mathrm{SL}_2(\R)$ acts linearly on $\mathcal{H}_n$ by
$$r_n(M)Q:=Q(ax+by,cx+dy),$$
and $r_n$ defines a representation $r_n:\mathrm{SL}_2(\R)\rightarrow \mathrm{SL}_{n+1}(\R)$. If $n$ is even then $r_n(-I)=r_n(I)$ and therefore
$r_n$ induces a natural representation of $\mathrm{PSL}_2(\R)$ into $\mathrm{SL}_{n+1}(\R)$. If $n$ is odd, we need to use a spin structure on $X=\Gamma\backslash \H$ given by a character
$\chi:\widetilde{\Gamma}\rightarrow \{\pm 1\}$ such that $\chi(-I)=-1$ and we define for all $\gamma \in \widetilde{\Gamma}$ such that $\pi(\widetilde{\gamma})=\gamma$ by
$$\widetilde{r_n}(\gamma):=(\chi (\widetilde{\gamma}))^n r_n(\widetilde{\gamma}).$$
Note that $\widetilde{r_n}$ now makes sense regardless of the parity of $n$ and defines a linear representation 
$$\widetilde{r_n}:\Gamma\rightarrow \mathrm{SL}_{n+1}(\R).$$
An easy calculation shows that if $\gamma \in \Gamma$, then we have
\begin{equation}
\label{traceirred}
\mathrm{tr}(\widetilde{r_n}(\gamma))= (\chi (\widetilde{\gamma}))^n\sum_{j=0}^n e^{(2j-n)\ell(\gamma)/2}.
\end{equation}
Using the fact that the matrix coefficients of $r_n(M)$ are homogeneous polynomials of degree $n$ in $a,b,c,d$, a simple compactness argument shows that for any
matrix norm $\Vert.\Vert$, one can find $C>0$ such that for all $\gamma \in \Gamma$,
\begin{equation}
\label{normeq}
\Vert \widetilde{r_n}(\gamma)\Vert \leq C \Vert \gamma \Vert^n,
\end{equation}
where $$\gamma \simeq \left ( \begin{array}{cc}a&b\\c&d   \end{array}    \right)$$
and
$$\Vert \gamma \Vert^2= a^2+b^2+c^2+d^2.$$
If $\rho:\Gamma\rightarrow \mathrm{PSL}_2(\R)$ is a Teichm\"uller representation, then clearly (\ref{normeq}) and (\ref{traceirred}) both show that $\rr=\widetilde{r_n}\circ \rho$ is a representation
of Teichm\"uller type with $\beta=n/2$. Notice that the trace formula (\ref{traceirred}) also shows that for $n=1,2$ these representations are equivalent to the ones previously defined.

By the same type of calculations as previously performed, one can deduce in the case $\rho=\Id$ that we have the identities:
$$Z_\Gamma(s,\rr)=\prod_{j=0}^n Z_\Gamma \left (s-j+\frac{n}{2}\right)\ \ \mathrm{if}\ n\ \mathrm{is\  even},$$
$$Z_\Gamma(s,\rr)=\prod_{j=0}^n Z_\Gamma \left (s-j+\frac{n}{2},\chi \right)\ \ \mathrm{if}\ n\ \mathrm{is\  odd}.$$
As before, this leads to explicit expressions for the spectrum of $\Delta_\rr$ in terms of Dirac or Laplace spectrum, depending on the parity of $n$, which perfectly generalizes to arbitrary dimensions
the results obtained before. It is possible to complicate the picture even more by choosing $p$ representations $\rr_1,\ldots,\rr_p$ such that
$$\rr_j=r_{d_j}\circ \rho:\Gamma \rightarrow \mathrm{SL}_{d_j}(\R),$$
and then define
$$\rr:=\rr_1\oplus \ldots \oplus \rr_p:\Gamma \rightarrow \mathrm{SL}_{d_1+\ldots d_p}(\R),$$
where $\rr_1\oplus \ldots \oplus \rr_p$ acts on $\C^{d_1}\times \ldots \C^{d_p}$ by 
$$\rr_1\oplus \ldots \oplus \rr_p(\gamma)(v_1,\ldots,v_p):=(\rr_1(\gamma)v_1,\ldots,\rr_p(\gamma)v_p).$$
When $\rho=\Id$, the spectrum of $\Delta_\rr$ can be computed as above and is composed of a mixture of Dirac and Laplace spectrum.

\begin{section}{Bowen-Series/Adler-Flatto coding and topological pressure}

\subsection{Boundary coding}
From now on, we will use the unit disc model $\mathbb{D}$, with boundary $S^1=\partial \mathbb{D}$ of the hyperbolic plane, endowed with the metric of constant negative curvature
$$ds^2=\frac{4dz\overline{dz}}{(1-\vert z\vert^2)^2}.$$ 
We will view the co-compact group $\Gamma$ as a subgroup of $\mathrm{PSU}(1,1)$ i.e., the set of Mobius transforms of the form
$$\gamma(z)=\frac{az+b}{\overline{b}z+\overline{a}},$$
where $\vert a\vert^2-\vert b\vert^2=1$. The action of $\Gamma$ on $S^1$ can be ``coded" by a single expanding, piecewise analytic map on $S^1$ often referred as a Bowen-Series
map associated to $\Gamma$. The construction is not unique and depends on a choice of a suitable fundamental domain with polygonal geodesic boundary and associated side-pairing isometries which generate $\Gamma$. See \cite{BS, Series, AdlerFlatto}, and also \cite{Pollicott1, PollRocha}.
More precisely, we have the following facts, which summarize the main points proved in the above cited papers.
\begin{propo}
\label{Coding} 
\begin{enumerate}
\item There exists a finite set of generators and their inverses $\mathcal{G}\subset \Gamma$.
 \item There exists a finite covering $S^1=\cup_{j=1}^k I_j$ by $k$ closed intervals \footnote{In the Adler and Flatto construction \cite{AdlerFlatto}, one can take $k=16g-8$, where $g$ is the genus of $X$.} $I_j$ such that $\mathrm{Int}(I_j)\cap \mathrm{Int}(I_i)=\emptyset$ for all $i\neq j$.
 \item There exists a piecewise analytic map $T:S^1\rightarrow S^1$ such that for all $j=1,\ldots,k$, $T\vert_{\mathrm{Int}(I_j)}=\gamma_j \in \mathcal{G}$ and $T$ is Markov:
 for each $I_j$, $T(I_j)$ is a finite union of intervals $I_i$.
 \item The map $T$ is eventually expanding i.e. , there exists $D>1$ such that $\inf_{x \in S^1}\vert (T^2)'(x) \vert \geq D$.
 \item  The map $T$ is topologically mixing: if $A$ is the $k\times k$ matrix defined by $A(i,j)=1$ if $T(I_i)\supset I_j$, $A(i,j)=0$ otherwise, then there exists $p_0>1$ such that
 $A^{p_0}$ has all its entries positive. 
\end{enumerate}
\end{propo}
We need more definitions. 
If $T^nx=x$, i.e. $x$ is a periodic point for $T$, we denote by $N(x)$ its {\it primitive period} i.e., the smallest integer $N\geq 1$ such that $T^Nx=x$. Obviously, $N(x)$ divides $n$.
Given an element $\gamma \in \Gamma \setminus \{\Id\}$, we will denote by $\vert \gamma\vert$ its word length with respect to the set of generators $\mathcal{G}$ which is 
$$\vert \gamma \vert=\min\{ N\ :\ \gamma=\gamma_1\ldots\gamma_N\ \mathrm{with}\ \gamma_j\in \mathcal{G}\}.$$
An admissible word $\alpha=(\alpha_1,\ldots,\alpha_n)\in \{1,\ldots,k\}^n$ is a word (of length $n=\vert \alpha\vert$) such that 
$$A(\alpha_i,\alpha_{i+1})=1\ \mathrm{for\  all}\  i=1,\ldots,n-1.$$ 
If $\alpha$ is an admissible word, we set 
$$g_\alpha:=\gamma_{\alpha_1}^{-1}\circ \ldots \circ \gamma^{-1}_{\alpha_n}.$$ 
For all $j$ such that $A(\alpha_n,j)=1$, $g_\alpha$ maps $I_j$ in
$I_{\alpha_1}$, and we have for all $x\in I_j$, $T^n(g_\alpha x)=x$ i.e. $g_\alpha$ is an inverse branch of $T^n$ with word length $\vert g_\alpha \vert=n$. 
The set of all inverse branches of $T^n$ corresponds to the set of mobius maps $\{ g_\alpha\ :\ \vert \alpha \vert=n,\ \alpha\ \mathrm{admissible}\}$. Each inverse branch $g_\alpha$ has a unique attracting fixed point denoted by $x_{g_\alpha}^-\in I_{\alpha_1}$.

\bigskip
The main interest in the above coding is the following property.
\begin{propo}
\label{primecoding}
With at most finitely many exceptions \footnote{The Boundary coding overcounts finitely many primitive geodesics, which corresponds to geodesics ending at boundary points of the $I_j$, see \cite{PollRocha, Pollicott1} for more comments.}, there is a bijection between:
\begin{itemize}
\item Primitive conjugacy classes $\{\gamma\}$ in $\Gamma$ with word length $n$ and translation length $\ell(\gamma)$.
\item Primitive periodic orbits $\{x,Tx,\ldots, T^{n-1}x\}$ with period $n$ and $\vert (T^n)'(x)\vert=e^{\ell(\gamma)}$. 
\item Primitive inverse branches $g_\alpha$ with $\vert \alpha \vert=n$ and $\vert g_\alpha'(x_{g_\alpha}^-)\vert=e^{-\ell(\gamma)}$.
\end{itemize}
\end{propo}

\subsection{Topological pressure and Mostow's map}
Le $I:=\sqcup_{i=1}^k I_i$ be the disjoint union of the intervals $I_i$. If $0<r<1$ we denote by $C^r(I)$ the space of piecewise $r$-H\"older functions on each $I_i$,
with respect to the usual angular distance $d_{S^1}(x,y)$ on $S^1$.
Given $\varphi \in C^r(I)$ for some $0<\alpha<1$, we denote by $P(\varphi)$ the {\it topological pressure}, defined by
$$P(\varphi)=\sup_{\mu} \left (  h_\mu(T)+\int_{S_1} \varphi d\mu     \right),$$
where the supremum runs over all $T$-invariant probability measures and $h_\mu(T)$ is the measure theoretic entropy of $\mu$. This supremum is realized by a unique probability measure $\mu_\varphi$ called the equilibrium
state of the potential $\varphi$, which is mixing, see \cite{PP} for example. Moreover, by \cite{PP} we have the following limit expression for the topological pressure:
$$e^{P(\varphi)}=\lim_{n\rightarrow +\infty} \left ( \sum_{T^n x=x} e^{\varphi^{(n)}(x)} \right)^{1/n},$$
where the sum runs over all $n$-periodic points of $T$ and $$\varphi^{(n)}(x)=\varphi(x)+\varphi(Tx)+\ldots+\varphi(T^{n-1}x).$$
In the latter, we will denote by $\tau(x):=\log\vert T'(x)\vert$ the distortion function for the Bowen-Series map $T$. By Proposition \ref{primecoding}, we know that
whenever $T^n x=x$ is a primitive periodic orbit, we have
$$\tau^{(n)}(x)=\log \vert (T^n)'(x)\vert=\ell(\gamma),$$
where $\gamma \in \mathcal{P}$ is in correspondence with the primitive periodic orbit $\{x,Tx,\ldots,T^{n-1}x\}$.
Alternatively, if $x_\alpha^-$ is the unique attracting fixed point of an inverse branch $\gamma=g_\alpha$, then we have
$$\log \vert (g_\alpha)'(x_\alpha^-)\vert=-\tau^{(n)}(x_\alpha^-)=-\ell(\gamma).$$
Given a Teichm\"uller representation $\rho:\Gamma\rightarrow \mathrm{PSU}(1,1)$, we would like to express the distorted translation lengths
$\ell_\rho(\gamma):=\ell(\rho(\gamma))$ via a suitable cocycle over the dynamical system $T:S^1\rightarrow S^1$. The main result of this section is as follows.
\begin{propo}
\label{Mostowmap} For each Teichm\"uller representation $\rho$, there exists $0<r\leq1$ and a H\"older function $\psi\in C^r(I)$ such that for all $\gamma \in \mathcal{P}$
with $\ell(\gamma)=\tau^{(n)}(x)$, with $T^n x=x$, we have
$$\ell_\rho(\gamma)=\psi^{(n)}(x).$$
Morevover, if $\rr:\Gamma\rightarrow \mathrm{GL}(V)$ is a Teichm\"uller like representation associated to $\rho$ with exponent $\beta$, the characteristic exponent $\delta(\rr)$
coincides with the unique zero of the map
$$\sigma \mapsto P(\beta \psi-\sigma \tau).$$
\end{propo}
Before we can give a proof, we need to recall a few facts. The Busemann cocycle $C_\xi(z,w)$ is defined (given $z,w \in \mathbb{D}$ and $\xi \in \partial \mathbb{D}$) by
$$C_\xi(z,w)=\lim_{t\rightarrow \infty} (d(z,\xi_t)-d(w,\xi_t)),$$
where $\xi_t$ is a parametrized geodesic line whose end point as $t\rightarrow \infty$ is $\xi \in S^1$. We summarize below the main useful properties of the Busemann cocycle in the disc model, valid for all $z,w,y \in \mathbb{D}$ and $\xi \in S^1$.
\begin{itemize}
\item $C_\xi(z,w)=-C_\xi(w,z)$.
\item $C_\xi(z,y)=C_\xi(z,w)+C_\xi(w,y)$.
\item $C_\xi(z,w)=\log \left ( \frac{(1-\vert w\vert^2)\vert z-\xi\vert^2}{(1-\vert z\vert^2)\vert w-\xi\vert^2}   \right).$
\item For all $\gamma \in \mathrm{PSU}(1,1)$, we have $C_\xi(\gamma z,\gamma w)=C_{\gamma^{-1}\xi}(z,w)$.
\item For all $\gamma \in \mathrm{PSU}(1,1)$, we have
$$e^{C_\xi(0,\gamma^{-1} 0)}=\frac{\vert \gamma(0)-\gamma(\xi)\vert^2}{1-\vert \gamma(0)\vert^2}=\vert \gamma'(\xi)\vert.$$
\end{itemize}
Basic Teichm\"uller theory, see \cite{Bers} pages 267-268, tells us that there exists a quasi-conformal homeomorphism $F_\rho:\mathbb{D}\rightarrow \mathbb{D}$ such that for all $\gamma \in \Gamma$, we have
$$\rho(\gamma)=F \circ \gamma \circ F^{-1}.$$
Moreover, $F_\rho$ extends as a homeomorphism of $\partial \mathbb{D}$ which is by Mori's theorem \cite{Mori} H\"older continuous. This boundary map $F$, often called Mostow's map, conjugates in a H\"older regular way the action of $\Gamma$ and $\rho(\Gamma)$ on the boundary $S^1$.
Notice that if $x_\gamma^+,x_\gamma^-$
are the repelling and attracting fixed points of $\gamma \in \Gamma$, then the repelling and attracting fixed points of $\rho(\gamma)$ are $F(x_\gamma^+)$ and $F(x_\gamma^-)$.
We then define $\psi$ as follows:
for all $x\in I_j$, set
$$\psi(x):=C_{F(x)}(0, \rho(\gamma_j)^{-1}0). $$
It is immediate to check that $\psi$ is indeed piecewise H\"older on $I$. Let $x\in S^1$ be such that
$T^nx=x$ and $x\in I_{i_1}$, $Tx\in I_{i_2}$,..., $T^{n-1}x \in I_{i_n}$ so that 
$$T^n x=\gamma_{i_n}\circ \ldots \circ \gamma_{i_1}x=\gamma(x),$$
with $x=x^+_\gamma$.
Writing
$$\psi^{(n)}(x)=\psi(x)+\ldots+\psi(T^{n-1}x)$$
$$= C_{F(x)}(0,\rho(\gamma_{i_1})^{-1}0)+C_{F(\gamma_{i_1}x)}(0,\rho(\gamma_{i_2})^{-1}0)+\ldots+
C_{F(\gamma_{i_{n-1}}\ldots \gamma_{i_1}x)}(0,\rho(\gamma_{i_n})^{-1}0),$$
and using the fact that $\rho(\gamma)\circ F=F\circ \gamma$, and the invariance property of the Busemann function by isometries, we get
$$\psi^{(n)}(x)=C_{F(x)}(0,\rho(\gamma_{i_1})^{-1}0)+C_{F(x)}(\rho(\gamma_{i_1})^{-1}0,\rho(\gamma_{i_2}\gamma_{i_1})^{-1}0)+\ldots$$
$$+C_{F(x)}(\rho(\gamma_{i_{n-1}}\ldots\gamma_{i_1})^{-1}0,\rho(\gamma_{i_{n}}\ldots\gamma_{i_1})^{-1}0).$$
The cocycle property then yields
$$\psi^{(n)}(x)=C_{F(x)}(0, \rho(\gamma_{i_{n}}\ldots\gamma_{i_1})^{-1}0)=\log\vert \rho(\gamma)'(F(x))\vert=\ell_\rho(\gamma),$$
and the first claim is proved. By the very definition of representations of Teichm\"uller type, we know that 
$$\delta(\rr)=\inf \left \{s>0\ :\ \sum_{\gamma \in \mathcal{P}, k\geq 1} e^{\beta k\ell_\rho(\gamma)-sk\ell(\gamma)}<\infty \right \}.$$
 We can use Proposition \ref{primecoding} and consider the sum
$$S_n(s):=\sum_{n\geq 1} \sum_{T^n x=x} \frac{1}{N(x)} e^{\beta \psi^{(n)}(x)-s\tau^{(n)}(x)}.$$ 
Using Fubini we can write formally
$$S_n(s)=\sum_{n\geq 1} \sum_{d\vert n}\frac{1}{d}\sum_{T^dx=x\atop N(x)=d}  e^{\frac{n}{d}(\beta \psi^{(d)}(x)-s\tau^{(d)}(x))}$$
$$=\sum_{k\geq 1} \sum_{d\geq 1} \frac{1}{d} \sum_{T^dx=x\atop N(x)=d}  e^{k(\beta \psi^{(d)}(x)-s\tau^{(d)}(x))}$$
$$=\sum_{k\geq 1} \sum_{\gamma \in \mathcal{P}} e^{\beta k \ell_\rho(\gamma)-sk\ell(\gamma)}+
\sum_{k\geq 1} \sum_{g\in \{g_1,\dots,g_m\} \subset \mathcal{P}} e^{\beta k \ell_\rho(g)-sk\ell(g)}.$$
where $g_1,\ldots,g_m$ are finitely many primitive conjugacy classes that are overcounted by the Bowen-series coding.
Clearly by the properties of the topological pressure, $S_n(s)$ converges if and only if $P(\beta \psi-s\tau)<0$.
In the above sum, the extra terms do not change the nature of the convergence, therefore $\delta(\rr)=\inf \{ s\ :\ P(\beta \psi-s\tau)<0\}$.
The topological pressure $\sigma \mapsto P(\beta \psi-\sigma\tau):=f(\sigma)$ is actually a real analytic function (\cite{PP}) which, as follows directly from the variational definition, is strictly decreasing and satisfies $\lim_{\sigma \rightarrow +\infty} f(\sigma)=-\infty$ while $\lim_{\sigma \rightarrow -\infty} f(\sigma)=+\infty$. As a conclusion, $f(\sigma)$ has a unique
zero which coincides with $\delta(\rr)$.
\end{section}

\begin{section}{Some technical bounds}
For all the admissible words $\alpha$ of length $n$ such that $A(\alpha_n,i)=1$, we denote by $\mathscr{C}_{\alpha\vee i}$ the ''cylinder set'' associated to the word $\alpha \vee i$
of length $n+1$ obtained by concatenation of $\alpha$ and $i$: 
$$\mathscr{C}_{\alpha\vee i}:=g_\alpha(I_i).$$
For all $n\geq 0$, cylinder sets form a "partition" of $S^1$:
$$S^1=\bigcup_{i=1}^k\bigcup_{\vert \alpha\vert=n,\atop A(\alpha,i)=1} \mathscr{C}_{\alpha\vee i},$$
and whenever $\alpha\vee i\neq \beta \vee j$, we have $\mathrm{Int}( \mathscr{C}_{\alpha\vee i})\cap \mathrm{Int}(\mathscr{C}_{\beta\vee j})=\emptyset $.

From the fact that the Bowen-Series map is eventually expanding, we know that there exist $C_0>0$ and $0<\theta_0\leq \theta_1<1$ such that for all $\alpha,i$ with
$A(\alpha,i)=1$, for all $x\in I_i$, we have
\begin{equation}
 \label{inversebranch1}
 C_0^{-1}\theta_0^n\leq \vert g_\alpha'(x)\vert \leq C_0 \theta_1^n.
\end{equation}
Given $\gamma \in \mathrm{PSU}(1,1)$, we can always choose a smooth lift $\widetilde{\gamma}:\R\rightarrow \R$ such that for all $\theta \in \R$, 
$\gamma(e^{i\theta})=e^{i\widetilde{\gamma}(\theta)}$. Using the fact that $\gamma$ is a Mobius map of the form
$$\gamma(z)=\frac{az+b}{\overline{a}z+\overline{b}}$$
with $\vert a\vert^2-\vert b\vert^2=1$,
we can deduce that 
\begin{equation}
\label{confderivative}
\forall\ \theta \in \R,\ \widetilde{\gamma}'(\theta)=\vert \gamma'(e^{i\theta}) \vert.
\end{equation}
Combining (\ref{inversebranch1}) and (\ref{confderivative}) we get that for all $\alpha,i$ with $A(\alpha,i)=1$ and $\vert \alpha\vert=n$, we have
\begin{equation}
 \label{Cylinderbound}
 \mathrm{diam}_{S^1}(\mathscr{C}_{\alpha\vee i})\geq \widetilde{C_0} \theta_0^n,
\end{equation}
for some uniform constant $\widetilde{C_0}>0$.
We are now ready to prove the following fact.
\begin{lem}(Separation)
\label{separation}
There exist $C_1>0$ and $0<\theta_2<1$ such that for all $n$ and all words $\alpha\neq \beta$ with $\vert \alpha\vert=\vert \beta\vert=n$, for all $i$ such that $A(\alpha_n,i)=A(\beta_n,i)=1$, we have for all $x\in I_i$ 
$$\vert g_\alpha(x)-g_\beta(x)\vert \geq C_1 \theta_2^n. $$
\end{lem}
\noindent {\it Proof.} First observe that if $\mathscr{C}_{\alpha \vee i} \cap \mathscr{C}_{\beta \vee i}=\emptyset$, then by (\ref{Cylinderbound}) we have 
$$\mathrm{dist}_{S^1}(\mathscr{C}_{\alpha \vee i},\mathscr{C}_{\beta \vee i})\geq  \widetilde{C_0} \theta_0^n,$$
and the claim is proved since both chordal and angular distance are equivalent on $S^1$. We can therefore restrict to the case when $\mathscr{C}_{\alpha \vee i} \cap \mathscr{C}_{\beta \vee i}\neq\emptyset$, in which case this intersection is a single point. Let us set
$$\mathscr{C}_{\alpha \vee i}=\{ e^{i\theta}\ :\ a\leq \theta \leq b\},\  \mathscr{C}_{\beta \vee i}=\{ e^{i\theta}\ :\ b\leq \theta \leq c\},$$
where $a<b<c$ and set also
$$I_i=\{ e^{i\theta}\ :\ \omega_1\leq \theta \leq \omega_2\},$$
with $\omega_1<\omega_2$.
We lift both maps $g_\alpha,\ g_\beta$ so that $g_\alpha(e^{i\theta})=e^{i\widetilde{g_\alpha}(\theta)}$, $g_\beta(e^{i\theta})=e^{i\widetilde{g_\beta}(\theta)}$, and
choose those lifts \footnote{This is possible  because these inverse branches are maps with degree $1$.} such that $\widetilde{g_\alpha}([\omega_1,\omega_2])=[a,b]$ while $\widetilde{g_\beta}([\omega_1,\omega_2])=[b,c]$.
For all $\omega_1\leq \theta \leq \omega_2$, we have
$$\widetilde{g_\beta}(\theta)-\widetilde{g_\alpha}(\theta)\geq \int_{\omega_1}^\theta \vert g_\beta'(\varphi)\vert d\varphi\geq (\theta-\omega_1)C_0^{-1}\theta_0^n,$$
while we also have
$$\widetilde{g_\beta}(\theta)-\widetilde{g_\alpha}(\theta)\geq (b-a)-\int_{\omega_1}^\theta  \vert g_\beta'(\varphi)\vert d\varphi\geq (b-a)-(\theta-\omega_1)C_0\theta_1^n.$$
Hence for all $\omega_1\leq \theta \leq \omega_2$ we get
$$\vert \widetilde{g_\beta}(\theta)-\widetilde{g_\alpha}(\theta)\vert \geq G_0(\theta):=\max \left\{(\theta-\omega_1)C_0^{-1}\theta_0^n, (b-a)-(\theta-\omega_1)C_0\theta_1^n \right\}.$$
A simple calculation then shows that for all $\omega_1\leq \theta \leq \omega_2$ and all $n$ large we have
$$G_0(\theta)\geq \frac{(b-a)C_0^{-1}\theta_0^n}{C_0^{-1}\theta_0^n+C_0\theta_1^n}\geq  \frac{\widetilde{C_0}C_0^{-1}\theta_0^{2n}}{C_0^{-1}\theta_0^n+C_0\theta_1^n}\geq C_1 \left ( \frac{\theta_0^2}{\theta_1}\right)^n,$$
and the proof is done. $\square$

We point out that we have used here in an essential way the fact that the Mobius maps $g_\alpha,\ g_\beta$ are orientation preserving and 
their derivatives enjoy exponential bounds. 

\bigskip
We need another technical fact. 
\begin{lem}(Norm Bound)
\label{nb1}
Let $\rr:\Gamma \rightarrow \mathrm{GL}(V)$ be a linear representation of Teichm\"uller type, with parameter $\beta$ and associated Teichm\"uller representation $\rho$.
Then for any matrix norm $\Vert .\Vert$ on $\mathrm{GL}(V)$, there exists $C_2>0$ such that for all $n$ and all admissible word $\alpha$ with length $n$ we have
$$\Vert \rr^{-1}( g_\alpha)\Vert\leq C_2 e^{\beta \ell_\rho( g_\alpha)}=C_2e^{\beta \psi^{(n)}(x_{g_\alpha}^+)}.$$
\end{lem}
\noindent {\it Proof.} Let $g_\alpha=\gamma_{\alpha_1}^{-1}\circ \ldots \circ \gamma_{\alpha_n}^{-1}$ be an inverse branch ot $T^n$ associated to an admissible word $\alpha$. 
Let $x_{g_\alpha}^-$ and $x_{g_\alpha}^+$ denote the attracting and repulsing fixed points of $g_\alpha$ which are distinct since $\Gamma$ is a co-compact group. 
By construction \footnote{See for example \cite{AdlerFlatto}, Thm 3.1, Thm 3.4 and Def 6.1. Both "left" and "right" Bowen-Series map work. See also \cite{PollRocha} for a picture of the side pairings and mapping properties.} 
of Bowen-Series maps \cite{BS, AdlerFlatto}, we have $x_{g_\alpha}^- \in I_{\alpha_1}$, 
while $x_{g_\alpha}^{+}\in I_j$ where $I_j\cap I_{\alpha_n}=\emptyset$. In particular $d_{S^1}(x_{g_\alpha}^- ,x_{g_\alpha}^+)\geq \min_{1\leq i \leq k} (\mathrm{diam}_{S^1}(I_i))$.
Since Mostow's map $F_\rho:S^1\rightarrow S^1$ is a continuous homomorphism, an immediate compactness argument shows that there exists $\eta>0$ such that for all $g_\alpha$
$$d_{S^1}(x_{\rho(g_\alpha)}^+,x_{\rho(g_\alpha)}^-)=d_{S^1}(F_\rho(x_{g_\alpha}^+),F_\rho(x_{g_\alpha)}^-)\geq \eta.$$
Given an hyperbolic isometry $\gamma$, in the disc model, the hyperbolic distance from the origin $0$ to the axis $A_\gamma$ of $\gamma$ is uniquely determined by the angular distance between the two fixed points $x_\gamma^+,\ x_\gamma^- \in S^1$. We therefore deduce that there exists $M>0$ such that for all $n$, for all inverse branch $g_\alpha$ we have 
$$d(A_{g_\alpha}, 0)\leq M.$$
Using property $1)$ of representations of Teichm\"uller type, we can write for all $g_\alpha$,
$$\Vert \rr^{-1}(g_\alpha)\Vert_1\leq Ce^{\beta d(\rho(g_\alpha^{-1}) 0,0)}=Ce^{\beta d(0, \rho(g_\alpha) 0)}.$$
By the triangle inequality, choosing a point $z\in A_{\rho(g_\alpha)}$, we have
$$d(0,\rho(g_\alpha)0)\leq  d(0,z)+d(z, \rho(g_\alpha) z)+d(\rho(g_\alpha) z, \rho(g_\alpha)0),$$
and thus
$$d(0,\rho(g_\alpha)0)\leq 2d(A_{\rho(g_\alpha)},0)+\ell_\rho(g_\alpha).$$
We can now conclude that
$$\Vert \rr^{-1}(g_\alpha)\Vert_1\leq Ce^{2\beta M+\beta \ell_\rho(g_\alpha)},$$
and the proof is done. $\square$
\end{section}

\begin{section}{Fredholm determinants and holomorphic function spaces} 

\subsection{Using zeta functions}
Following Ruelle \cite{Ruelle1} and Pollicott \cite{Pollicott1}, since each isometry $\gamma_i:I_i\rightarrow S^1$ is real-analytic, and because of the contraction properties, one can find for each $i$ an open set $\C\supset U_i\supset I_i$
such that whenever $A(i,j)=1$,
$$\overline{\gamma_i^{-1}(U_j)}\subset U_i.$$
Here $\gamma_i^{-1}:U_j\rightarrow U_i$ is simply understood as the Mobius map $\gamma_j^{-1}$ (which is holomorphic on a neighborhhod of $S^1$) acting on the complex open set $U_j$.
Since the mappings $\gamma_i^{-1}:U_j\rightarrow U_i$ are {\it holomorphic contractions}, it is not surprising that they have strong spectral properties when acting on ``favorable function spaces".  
For a precise definition of ``favorable function spaces"  and possible choices we recommend to read \cite{BJ1}.

Recall that $V$ denotes the representation space
of $\rr:\Gamma \rightarrow \mathrm{GL(V)}$, which is endowed with a Hermitian metric $\langle .,.\rangle_V$, with associated norm $\Vert . \Vert_V$.
 In our case, we will use Bergman spaces $H^2(U_i)$ defined as
$$H^2(U_i,V):=\left\{ f:U_i\rightarrow V\ :\ f\ \mathrm{is\ holomorphic\ and}\ \int_{U_i}\Vert f(z)\Vert^2_Vdm(z)<\infty \right \},$$
where $m$ denotes Lebesgue measure on $\C$. Each space $H^2(U_i,V)$ is a Hilbert space endowed with the $L^2$-norm.
We now set
$$\mathcal{H}^2:=\bigoplus_{i=1}^k H^2(U_i,V)$$
And define a transfer operator $\lt_{s,\rr}: \mathcal{H}^2\rightarrow \mathcal{H}^2$ by 
$$\forall\ z\in U_j,\ \lt_{s,\rr}(F)(z):=\sum_{A(i,j)=1} [(\gamma_i^{-1})']^s(z) \rr^{-1}(\gamma_i^{-1})F\circ \gamma_i^{-1}(z),$$
where  $[(\gamma_i^{-1})'](z)$ is defined by analytic continuation to $U_j$ of $\vert (\gamma_i^{-1})'\vert$, defined on $I_j$. The complex power $[(\gamma_i^{-1})']^s(z)$ is then defined by taking a suitable holomorphic branch of the logarithm such that the restriction to $I_j$ coincides with the usual complex power of a real number, see the next subsection for details.
The $n$-th iterate $\lt^n_{s,\rr}$ acts as follows:
$$\forall\ z\in U_j,\ \lt_{s,\rr}^n(F)(z):=\sum_{A(\alpha_{n},j)=1\ \atop \vert \alpha\vert=n}[g_\alpha']^s(z) \rr^{-1}(g_\alpha)F\circ g_\alpha(z),$$
where the sum runs over all admissible words of length $n$. By the general facts on Holomorphic contraction systems \cite{BJ1}, we know that $\lt_{s,\rr}$ is a compact trace class
operator when acting on $\mathcal{H}^2$. In addition, the trace is given by
$$\mathrm{tr}(\lt^n_{s,\rr})=\sum_{j=1}^k \sum_{A(\alpha_{n},j)=1\ \atop \vert \alpha\vert=n} 
\mathrm{tr}(\rr^{-1}(g_\alpha))\frac{\vert g_\alpha'\vert^s(x_{g_\alpha}^-)}{1-g_\alpha'(x_{g_\alpha}^-)}.$$

Lefschetz type formulas in this context date back to Ruelle \cite{Ruelle1}, see also \cite{Pollicott1}. For a proof in the vector valued case, including non-unitary representations, we refer the reader to \cite{PF}. A direct consequence of this trace formula combined with Proposition \ref{primecoding} is the following key fact: for all $s\in \C$, we have
$$\det(I-\lt_{s,\rr})= Z_\Gamma(s,\rr) \prod_{k\geq 0} \prod_{g_1,\ldots,g_m} \det(I-\rr(g_j)e^{-(s+k)\ell(g_j)}),$$
where $\{g_1,\ldots,g_m\}\subset \mathcal{P}$ is a finite set of primitive conjugacy classes that are "overcounted" by the Bowen-Series coding, see \cite{PollRocha,Pollicott1}.
In particular, zeros of $Z_\Gamma(s,\rr)$ are a subset of zeros of the Fredholm determinant $\det(I-\lt_{s,\rr})$. Notice also that for all $n\geq 1$, zeros of $\det(I-\lt_{s,\rr})$ are themselves, with multiplicity, a subset of the zero set of $\det(I-\lt_{s,\rr}^n)$, by using the product formula for trace class operators:
$$\det(I-\lt_{s,\rr}^n)=\det(I-\lt_{s,\rr})\det(I+\lt_{s,\rr}+\ldots+\lt_{s,\rr}^{n-1})$$
(see \cite{Simon}).
The main result of this paper (improved Weyl law) will follow from the next statement.
\begin{thm}
\label{main3}
 Let $\delta_0$ be the unique zero of the map $\sigma \mapsto P(2\beta \psi-2\sigma \tau)$, where the weights $\psi$ and $\tau$ are defined as in $\S 4$. For all $\sigma_0>\delta_0$, there exist constants $\kappa>0$, $M_0>0$
 and $\eta_0, T_0>0$ such that for $n(T)=2[\kappa \log T ]$, for all $s=\sigma+it$ with $\sigma\geq \sigma_0$ and $T_0\leq \vert t\vert\leq T$, 
 $$\log \vert \det(I -\lt_{s,\rr}^{n(T)})\vert \leq M_0 T^{1-\eta_0}.$$
 \end{thm}
First, we explain why we have $\delta_0<\delta$, the fact that $\delta_0>\delta/2$ will follow from the bounds proved in the last subsection (see geodesic stretch and bounds on $\delta_0$). By the theory of equilibrium states for Markov expanding maps (or subshifts of finite type), we know
that
$$P(2\beta \psi-2\delta_0 \tau)=0=h_{\mu_0}(T) +2\beta \int \psi d\mu_0-2\delta_0\int \tau d\mu_0,$$
where $\mu_0$ is a $T$-invariant probability measure associated to the H\"older potential $2\beta \psi-2\delta_0 \tau$. This equilibrium state $\mu_0$ has positive entropy
$h_{\mu_0}(T)>0$ by a general fact on equilibrium states, see \cite{Vc1}. Therefore, we have
$$0<2\left(h_{\mu_0}(T)+\beta \int \psi d\mu_0-\delta_0\int \tau d\mu_0 \right)$$
$$\leq 2 \sup_{\mu} \left (   h_{\mu}(T)+\beta \int \psi d\mu-\delta_0\int \tau d\mu  \right)=2P(\beta \psi-\delta_0 \tau).$$
Hence we get $P(\beta \psi-\delta_0 \tau)>0$ which implies that $\delta_0< \delta$. It is not difficult to see that $\delta_0$ is positive. For more refined explicit bounds on $\delta_0$ involving $\delta$ and the ''geodesic stretch", see the end of this section.

\bigskip 
Let us now show how this bound implies Theorem \ref{main2}. First observe that if $\lambda=s(1-s)$ with $s=\sigma+it$ satisfies 
$\lambda \in \mathcal{C}_\delta \setminus \mathcal{C}_{\sigma_0}$ then we can assume that
$$\sigma_0\leq \sigma \leq \delta\ .$$
Using the fact that asymptotically as $\vert t\vert \rightarrow \infty$ we have $\sqrt{\vert \lambda\vert}=\vert t\vert+O(1)$, we deduce that 
$$\#\{ \lambda \in \mathrm{Sp}(\Delta_{\rr})\ :\ \lambda \in \mathcal{C}_{\delta}\setminus \mathcal{C}_{\sigma_0}\ \mathrm{and}\ r\leq \vert \lambda \vert \leq r+\sqrt{r} \}$$
$$\leq \#\{ Z_\Gamma(s,\rr)=0\ :\ \sigma_0\leq \sigma \leq \delta\ \mathrm{and}\ \sqrt{r}-M\leq \vert t\vert \leq \sqrt{r}+M\},$$
for some constant $M>0$. To bound from above the number of zeros of $Z_\Gamma(s,\rr)$ in this region as $r\rightarrow +\infty$, we use Jensen's formula as below.
\begin{lem}
\label{Jensen}
 Let $f$ be a holomorphic function on the open disc $D(w,R)$, and assume that $f(w)\neq 0$. let $N_f(r)$ denote the number of zeros of $f$ in the closed disc $\overline{D}(w,r)$. For all $\widetilde{r}<r<R$, we have 
$$N_f(\widetilde{r})\leq \frac{1}{\log(r/\widetilde{r})} \left (  \frac{1}{2\pi} \int_0^{2\pi} \log \vert f(w+re^{i\theta})\vert d\theta-\log\vert f(w)\vert   \right).$$
\end{lem}
For a reference, one can look for example at \cite{Ransford}, Corollary 4.5.2.
Let us set
$$\mathcal{R}_0^{\pm}:= \{ s=\sigma+it\in \C\  :\ \sigma_0\leq \sigma \leq \delta\ \mathrm{and}\ \pm \sqrt{r}-M\leq  \pm t  \leq \pm \sqrt{r}+M\}.$$
We treat the case of $\mathcal{R}_0^+$, the case of $\mathcal{R}_0^-$ being similar.
Choose $\epsilon>0$ so that $\delta_0< \sigma_0-\epsilon<\sigma_0$. Set $R=\vert A-\sigma_0-iM\vert+\epsilon/3$ and $R'=\vert A-\sigma_0-iM\vert+2\epsilon/3$.
By fixing $A$ large enough, we can ensure that for all $r$ large enough we have
$$\mathcal{R}_0^+\subset D(A+i\sqrt{r},R)\subset D(A+i\sqrt{r},R')\subset \{ \Re(s)\geq \sigma_0-\epsilon\}.$$
In order to apply Jensen's formula to $\det(I -\lt_{s,\rr}^{n(T)})$ with $T=\sqrt{r}+R'$, we need to control from below
$$\log\vert \det(I -\lt_{A+i\sqrt{r},\rr}^{n(T)}) \vert.$$
By the Lefschetz trace formula, we have for large $\Re(s)$
$$\det(I -\lt_{s,\rr}^{n(T)})=\exp\left (-\sum_{p=1}^\infty \frac{1}{p} \mathrm{tr}( \lt_{s,\rr}^{pn(T)})\right),$$
with 
$$\mathrm{tr}( \lt_{s,\rr}^{pn(T)})= \sum_{j=1}^k \sum_{A(\alpha_{\vert \alpha \vert},j)=1\ \atop \vert \alpha\vert=pn(T)} 
\mathrm{tr}(\rr^{-1}(g_\alpha))\frac{\vert g_\alpha \vert'^s(x_{g_\alpha}^-)}{1-g_\alpha'(x_{g_\alpha}^-)}.$$
Using property $(2)$ from representations of Teichm\"uller type, we have
$$\vert \mathrm{tr}( \lt_{s,\rr}^{pn(T)})\vert \leq C \sum_{T^{pn(T)}x=x} e^{\beta \psi^{(pn(T))}(x)-\Re(s)\tau^{(pn(T))}(x)}.$$
In particular, if $\Re(s)=A$ is large enough, then the topological pressure $P(\beta \psi-A\tau)$ is negative and we get a positive lower bound on $\vert \det(I -\lt_{s,\rr}^{n(T)})\vert $
which is uniform in $r$. Applying the above Lemma on the disc $D(A+i\sqrt{r},R')$ to the holomorphic function $\det(I -\lt_{s,\rr}^{n(T)})$ we get
$$\# (\{ Z_\Gamma(s,\rr)=0\}\cap \mathcal{R}_0^+)=O(T^{1-\eta_0})+O(1), $$
and thus 
$$\#\{ \lambda \in \mathrm{Sp}(\Delta_{\rr})\ :\ \lambda \in \mathcal{C}_{\delta}\setminus \mathcal{C}_{\sigma_0}\ \mathrm{and}\ r\leq \vert \lambda \vert \leq r+\sqrt{r} \}=
O(T^{1-\eta_0})$$
$$=O(r^{1/2-\eta_0/2}),$$
which is the bound of Theorem \ref{main2}.

\subsection{Refined function space and proof of Theorem \ref{main3}}
We will need in the sequel to refine the function space $\mathcal{H}^2$. Let $h>0$ denote a small parameter. For each $j=1,\ldots,k$ set
$$U_j(h):=I_j+D(0,h),$$
where $D(0,h)$ is the open euclidean disc of radius $h$ centered at $0$ in $\C$. Then $U_j(h)$ is a simply connected open subset of $\C$. It is clear thar for all $h>0$ small enough we have
for all $j=1,\ldots,k$, $U_j(h)\subset U_j$. We then set
$$\mathcal{H}^2_h:= \bigoplus_{i=1}^k H^2(U_i(h),V),$$
endowed with the $L^2$ norm. These $h$-dependent families of Bergman spaces have been first used in the work of Guillop\'e-Lin-Zworski \cite{GLZ} for resonances of Schottky manifolds, where they considered families
of complex balls of radius $h$. In our analysis we will instead have to consider tubular neigborhoods of intervals. 
\begin{propo}
\label{morebounds}
\begin{enumerate}
 \item There exist $C>0$ and $0<\theta_3<1$ such that for all $i=1,\ldots,k$, for all admissible word $\alpha$ with $\vert \alpha\vert=n$ and $A(\alpha_n,i)=1$, we have
 $$\sup_{z\in U_i}\vert g'_\alpha(z)\vert\leq C_3\theta_3^n. $$
 \item There exists $n_0\geq 1$ such that for all $h>0$ small enough and all $n\geq n_0$, $\lt_{s,\rr}^n$ is a well defined compact trace class operator on $\mathcal{H}^2_h$. 
 In addition, if $\vert \alpha \vert=n\geq n_0$ and $g_\alpha$ maps $U_i(h)$ into $U_{\alpha_1}(h)$, we have for all $z\in U_i(h)$, 
 $$\mathrm{dist}(g_\alpha(z),\partial U_{\alpha_1}(h))\geq \frac{h}{2}. $$
 Moreover, the traces (hence the Fredholm determinant) do not depend on $h$. 
 \item There exists $B_0>0$ such that for all admissible word $\alpha$ with $\vert \alpha\vert=n$ and indices $i,j$ with $A(\alpha_n,i)=A(\alpha_n,j)=1$, we have for all $z\in U_i$ and $w\in U_j$
 $$ \left \vert  \frac{g'_\alpha(z)}{g'_\alpha(w)} \right \vert \leq B_0.$$
 \item There exists $B_1(\sigma_0)>0$ such that 
 for all $0\leq \Re(s)\leq M_0$, for all $i=1,\ldots,k$, for all admissible word $\alpha$ with $\vert \alpha\vert=n$ and $A(\alpha_n,i)=1$, we have for $h\leq \vert \Im(s)\vert^{-1}$ and
 $z\in U_i(h)$, 
 $$\left \vert  [g'_\alpha]^s(z)   \right \vert \leq B_1 \vert g'_\alpha(z)\vert^{\Re(s)}.$$
\end{enumerate}
\end{propo}
\noindent {\it Proof}. These estimates are rather standard for expanding Markov maps. Nevertheless we give some details. 
\begin{enumerate}
\item We can view $\Gamma$ as a discrete subgroup of $\mathrm{PSL}_2(\C)$. Each $\gamma \in \Gamma$ has an isometric circle $C_\gamma:=\{ z\in \C\ :\ \vert \gamma'(z) \vert=1\}$, with radius $r_\gamma$ and center $z_\gamma=\gamma^{-1}(\infty)$. 
Since $\Gamma$ is Fuchsian and discrete, the radii $r_\gamma$ can only accumulate at $0$ and the centers $x_\gamma$ can only accumulate on $S^1$. Notice that $x^+_\gamma \in D_\gamma:=\{ z\in \C\ :\ \vert \gamma'(z) \vert\geq 1\}$. As a consequence, there exists $n_0$ such that for all admissible word $\alpha$ with $n=\vert \alpha\vert\geq n_0$, then
$$D_{g_\alpha}\cap \bigcup_{A(\alpha_n,j)=1} U_j=\emptyset.$$
In particular, if $\vert \alpha\vert=n_0$, 
$$\sup_{z\in \cup_{A(\alpha_n,j)=1} U_j } \vert g'_\alpha(z)\vert<1.$$
The claim now follows from a simple iteration argument.
\item Let $U_i(h)$ be defined as above and pick an inverse branch $g_\alpha$ with word length $n$, such that $A(\alpha_n,i)=1$. Let $z\in U_i(h)$, and choose 
$\widetilde{z}\in I_i$ such that $\vert z-\widetilde{z}\vert<h$. By $1)$, we have
$$\vert g_\alpha(z)-g_\alpha(\widetilde{z})\vert \leq C\theta_3^nh,$$
and since $g_\alpha(\widetilde{z})\in I_{\alpha_1}$, we get that $g_\alpha(z) \in U_{\alpha_1}(h)$ and 
$$\mathrm{dist}(g_\alpha(z),\partial U_{\alpha_1}(h))\geq h/2,$$ as long as $C\theta_3^n\leq 1/2$. Therefore the operator $\lt_{s,\rr}^n$ is well defined on $\mathcal{H}_h^2$ for all $n\geq n_0$. Because it acts by uniform contractions on holomorphic functions, it is a trace class operator. The fact that traces do not depend on $h$ is because they are expressed only in term of fixed points which belong to the intervals $I_j$'s. By analytic continuation, the determinant is also independent of $h$.
\item Follows directly from $1)$ after taking logarithms and applying the chain rule. 
\item This is similar to \cite{Naud1, JN1}, however since we work with the disc model, there are some slight differences that we need to address. Given
$\gamma \in \Gamma$ with 
$$\gamma(z)=\frac{az+b}{\overline{a}z+\overline{b}},$$
we have for all $z\in S^1$,
$$\vert \gamma'(z)\vert=\frac{1}{(\overline{a}z+\overline{b})(az^{-1}+b)},$$
so that the analytic continuation $[\gamma'](z)$ is given by the same formula (which is holomorphic in a neighborhood of $S^1$).
Given an inverse branch $g_\alpha$ with $\vert \alpha \vert=n$ and $A(\alpha_n,i)=1$, $z=re^{i\theta}\in U_i(h)$, we write
$$[g_\alpha'](re^{i\theta})=\frac{1}{\vert \overline{a}e^{i\theta}+\overline{b}\vert^2+E(a,b,r,\theta)},$$
where
$$ E(a,b,r,\theta)=\overline{b}ae^{-i\theta}\left(\frac{1}{r}-1\right)+b\overline{a}e^{i\theta}(r-1).$$
It is immediate to check that for all $h$ small enough we have
$$\vert E(a,b,r,\theta)\vert \leq 3\vert ab\vert h.$$
Observe now that 
$$\frac{\vert E(a,b,r,\theta)\vert}{\vert \overline{a}e^{i\theta}+\overline{b}\vert^2}\leq 3h\frac{\vert \frac{b}{a}\vert}{\vert e^{i\theta}+\frac{\overline{b}}{\overline{a}}\vert^2}.$$
Since $-\frac{\overline{b}}{\overline{a}}$ is the center of the isometric circle of $g_\alpha$ and $z=re^{i\theta}\in U_i(h)$, this quantity is bounded uniformly as
\begin{equation}
\label{keybound1}
\frac{\vert E(a,b,r,\theta)\vert}{\vert \overline{a}e^{i\theta}+\overline{b}\vert^2}\leq \widetilde{C}h.
\end{equation}
In particular, from (\ref{keybound1}), we get that for all $h$ small enough, for all $z\in U_i(h)$
$$\vert [g'_\alpha](z)\vert\leq  A_1 \vert g'_\alpha(z)\vert,$$
where $A_1>0$ is some uniform constant. Going back to $[g'_\alpha]^s(z)$, we have
$$[g'_\alpha]^s(z)=e^{s \mathbb{L}([g'_\alpha](z))},$$
where $\mathbb{L}(w)$ is the principal holomorphic determination of the logarithm defined on $\C\setminus \R^-$ given by
$$\mathbb{L}(w)=\log\vert w\vert+i\mathrm{Arg}(w),$$
where $\mathrm{Arg}(w)\in (-\pi,+\pi)$ is the principal argument of $w$.
Writing
$$\vert [g'_\alpha]^s(z) \vert =e^{\Re(s) \log\vert [g'_\alpha](z)\vert}e^{-\Im(s) \mathrm{Arg}([g'_\alpha](z))},$$
we use (\ref{keybound1}) to bound the argument as
$$\vert \mathrm{Arg}([g'_\alpha](z))\vert \leq A_2 h.$$
It is now clear that if $\vert \Im(s)\vert\leq h^{-1}$ and $0\leq \Re(s)\leq M_0$, then we get
$$\vert [g'_\alpha]^s(z) \vert\leq A_3(M_0) \vert g_\alpha'(z)\vert^{\Re(s)},$$
for some $A_3(M_0)>0$ and the proof is done. $\square$
\end{enumerate}
Before we can give a proof of Theorem \ref{main3}, we need to recall a few facts on Bergman kernels which are essential at the technical level. If $\Omega$ is a bounded open set in $\C$, then the Hilbert space $H^2(\Omega,\C)$ has a reproducing kernel i.e., for all $f\in H^2(\Omega,\C)$, we have
$$f(z)=\int_\Omega B_\Omega(z,w)f(w)dm(w),$$
where $B_\Omega(z,w)$ is a smooth kernel on $\Omega\times \Omega$, called the Bergman kernel of $\Omega$. Given a Hilbert basis $(\phi_l)_{l\in \N}$ of 
$H^2(\Omega,\C)$, we have for all $z,w\in \Omega$,
$$\sum_{l\in \N}\phi_l(z)\overline{\phi_l(w)}=B_\Omega(z,w),$$
the convergence being uniform on every compact subset of $\Omega \times \Omega$, see for example Krantz \cite[Chapter 1]{Krantz}. In the latter, we will work on narrow domains
$\Omega=U_i(h)$ with $h\rightarrow 0$. We will need to use the following bounds on kernels $B_{U_i(h)}(z,w)$ which are proved in the Appendix (Proposition A.1).
\begin{itemize}
\item There exists a constant $M>0$ independent of $h$ such that for all $z,w\in U_i(h)$ with $\mathrm{dist}(z,\partial U_i(h))\geq h/2$ and $\mathrm{dist}(w,\partial U_i(h))\geq h/2$,
we have for all $h>0$ small enough $$\vert B_{U_i(h)}(z,w)\vert \leq M h^{-2}.$$
\item If in addition we have $\vert z-w\vert\geq ch^{\nu}$ for some $c>0$ and $0<\nu<1$, then we have for all $h>0$ small enough
$$\vert B_{U_i(h)}(z,w)\vert \leq M h^{-1-\nu}.$$
\end{itemize}
We are now ready to prove Theorem \ref{main3}. We assume that $s=\sigma+it$ with $\sigma \geq \sigma_0>\delta_0$ and $\vert t \vert \leq T$. We set 
$n(T)=2[\kappa\log(T)]$ where $T$ is large and $\kappa>0$ will be adjusted later on. We will let $\lt_{s,\rr}^{n(T)}$ act on the function space $\mathcal{H}^2_h$ with $h=T^{-1}$.
By a standard bound on Fredholm determinants, see \cite{Simon}, we know that 
$$\log \vert \det(I-\lt_{s,\rr}^{n(T)})\vert \leq \left \Vert \lt_{s,\rr}^{n(T)/2} \right \Vert_{HS}^2,$$
where $\Vert.\Vert_{HS}$ denotes the Hilbert-Schmidt norm, when acting on the function space $\mathcal{H}_h^2$. We now need to compute and estimate this Hilbert-Schmidt norm.
For each $i=1,\ldots,k$, let $(\phi_l^{(i)})_{l\in \N}$ be a Hilbert basis of $U_i(h)$. Let $\eb_1,\ldots,\eb_d$ be an orthonormal basis of the representation space $V_\rr$ which is endowed with a hermitian inner product $\langle .,. \rangle_V$. Then the family
$$\left ( \phi^{(i)}_l(z) \eb_m\right)_{l\in \N,\atop {1\leq i \leq k,\atop  1\leq m \leq d}},$$
is a Hilbert basis of $\mathcal{H}_h^2=\bigoplus_i H^2(U_i(h),V)$ and we have
$$\left \Vert \lt_{s,\rr}^{n(T)/2} \right \Vert_{HS}^2=\sum_{i,m}\sum_l \Vert \lt_{s,\rr}^{n(T)/2}(\phi_l^{(i)}\eb_m)\Vert_{\mathcal{H}^2_h}^2,$$
where $\phi_l^{(i)}\eb_m \in \mathcal{H}^2_h$ is understood as being
$$\phi_l^{(i)}(z)\eb_m =\left \{ \phi_l^{(i)}(z)\eb_m\ \mathrm{if}\ z\in U_i(h) \atop 0\ \mathrm{otherwise.} \right.$$
Setting for conveniency $\widetilde{n}(T)=n(T)/2$, we therefore have
$$\left \Vert \lt_{s,\rr}^{\widetilde{n}(T)} \right \Vert_{HS}^2=$$
$$\sum_{l,m,i,j} \int_{U_j(h)} \sum_{A(\alpha_{\widetilde{n}},j)=1 \atop A(\beta_{\widetilde{n}},j)=1} [g'_\alpha]^s(z)  \overline{[g'_\beta]^s(z)} 
\phi_l^{(i)}(g_\alpha z)\overline{\phi_l^{(i)}(g_\beta z)} \langle\rr^{-1}(g_\alpha)\eb_m, \rr^{-1}(g_\beta)\eb_m \rangle_V dm(z).$$
Using Fubini (justified by uniform convergence on compact sets for the infinite sum over $l$) we get
$$ \left \Vert \lt_{s,\rr}^{\widetilde{n}(T)} \right \Vert_{HS}^2=$$
$$\sum_j \int_{U_j(h)} \sum_{A(\alpha_{\widetilde{n}},j)=1 \atop A(\beta_{\widetilde{n}},j)=1}  \mathrm{tr}(\rr^{-1}(g_\alpha)\rr^{-1}(g_\beta)^*)[g'_\alpha]^s(z)  \overline{[g'_\beta]^s(z)} 
B_{U(h)}(g_\alpha z,g_\beta z)dm(z),$$
Where $U(h)=\sqcup_{i=1}^kU_i(h)$ and $B_{U(h)}(z,w)$ is the reproducing kernel of $\mathcal{H}^2_h$.
We then use the fact that (Cauchy-Schwartz inequality)
$$\vert \mathrm{tr}(\rr^{-1}(g_\alpha)\rr^{-1}(g_\beta)^*)\vert \leq \Vert \rr^{-1}(g_\alpha)\Vert_{HS} \Vert \rr^{-1}(g_\beta)\Vert_{HS} ,$$
where $\Vert.\Vert_{HS}$ is the Hilbert-Schmidt norm on $\mathrm{GL}(V)$ associated to the inner product on $V$. Using the norm bound from Lemma \ref{nb1} combined with
Proposition \ref{morebounds} (4), we get
$$\left \Vert \lt_{s,\rr}^{\widetilde{n}(T)} \right \Vert_{HS}^2\leq C \times $$
$$\sum_j\int_{U_j(h)} \sum_{A(\alpha_{\widetilde{n}},j)=1 \atop A(\beta_{\widetilde{n}},j)=1}  e^{\beta (\ell_\rho(g_\alpha)+\ell_\rho(g_\beta))} \vert g'_\alpha(z)\vert^{\Re(s)}\vert g'_\beta(z)\vert^{\Re(s)}
\vert B_{U(h)}(g_\alpha z,g_\beta z)\vert dm(z).$$
We will now split the above sum into two parts $\mathcal{S}_1(h)+\mathcal{S}_2(h)$, where
$$\mathcal{S}_1(h)= \sum_j\int_{U_j(h)} \sum_{A(\alpha_{\widetilde{n}},j)=1}  e^{2\beta \ell_\rho(g_\alpha)} \vert g'_\alpha(z)\vert^{2\Re(s)}
\vert B_{U(h)}(g_\alpha z,g_\alpha z)\vert dm(z),$$
which we call the ``diagonal part'' and 
$$\mathcal{S}_2(h)=\sum_j\int_{U_j(h)} \sum_{\alpha \neq \beta}  e^{\beta (\ell_\rho(g_\alpha)+\ell_\rho(g_\beta))} \vert g'_\alpha(z)\vert^{\Re(s)}\vert g'_\beta(z)\vert^{\Re(s)}
\vert B_{U(h)}(g_\alpha z,g_\beta z)\vert dm(z),$$
which we call the ``off-diagonal'' part. We first deal with the diagonal part. Using Proposition \ref{morebounds} (3), we have
$$\mathcal{S}_1(h)\leq C \sum_j\int_{U_j(h)} \sum_{A(\alpha_n,j)=1}  e^{2\beta \ell_\rho(g_\alpha)-2\sigma \ell(g_\alpha)} 
\vert B_{U(h)}(g_\alpha z,g_\alpha z)\vert dm(z).$$
Using the crude bound on Bergman kernels combined with Proposition \ref{morebounds} (2), we get uniformly in $z$
$$ \vert B_{U(h)}(g_\alpha z,g_\alpha z)\vert\leq C h^{-2},$$
and since we have $m(U(h))\leq Ch$ (either by direct computation or by using a general result on volumes of tubular neighborhood of smooth curves) we obtain that for $\Re(s)=\sigma\geq \sigma_0$, 
$$\mathcal{S}_1(h)\leq C h^{-1}\sum_j  \sum_{A(\alpha_{\widetilde{n}},j)=1}  e^{2\beta \ell_\rho(g_\alpha)-2\sigma_0 \ell(g_\alpha)},$$
where the words are of length $\widetilde{n}(T)=[\kappa \log T]$. Fixing $\epsilon >0$ small enough, this term can be controlled via the topological pressure for all $\widetilde{n}(T)\geq n_0$:
$$\mathcal{S}_1(h)\leq C Te^{\widetilde{n}(T)(P(\sigma_0)+\epsilon)},$$
where $P(\sigma_0)$ is a shorthand for $P(2\beta \psi-2\sigma_0\tau)$. Clearly, as long as $\sigma_0>\delta_0$ we have $P(\sigma_0)<0$ and we get a polynomial gain of magnitude
$O(T^{-\eta_0})$ for some $\eta_0>0$. 

We now deal with the ``off-diagonal sum'', and this is where the ``off-diagonal estimate'' on Bergman kernels is critical.
Observe that if $z\in U_j(h)$ and $\alpha,\beta$ are two distinct admissible words of length $\widetilde{n}(T)$, then picking $\widetilde{z}\in I_j$ with $\vert z-\widetilde{z}\vert <h$, we have by Lemma \ref{separation} and Proposition \ref{morebounds} (1),
$$\vert g_\alpha(z)-g_\beta(z)\vert \geq C_1\theta_2^{\widetilde{n}(T)}-2 C_3\theta_3^{\widetilde{n}(T)}h.$$
Taking $\kappa>0$ small enough so that $\theta_2 ^{\widetilde{n}(T)}\geq \sqrt{h}$, we end up with 
$$ \vert g_\alpha(z)-g_\beta(z)\vert \geq C_4 \sqrt{h},$$
for all $z\in U_j(h)$ and all $h=T^{-1}$ small enough. Using the off-diagonal bound on Bergman kernels and Proposition \ref{morebounds} (3), we get
$$ \mathcal{S}_2(h)\leq Ch^{-1-1/2}\sum_j\int_{U_j(h)} \sum_{\alpha \neq \beta}  e^{\beta (\ell_\rho(g_\alpha)+\ell_\rho(g_\beta))-\sigma_0\ell(g_\alpha)-\sigma_0\ell(g_\beta)}
 dm(z).$$
Using the topological pressure formula, we get (again fixing $\epsilon>0$ and taking $T$ large enough)
 $$\mathcal{S}_2(h)\leq CT^{1/2} e^{(2P(\beta \psi-\sigma_0 \tau)+\epsilon)\widetilde{n}(T)}.$$
 We then choose again $\kappa>0$ small enough so that 
 $$e^{(2P(\beta \psi-\sigma_0 \tau)+\epsilon)\widetilde{n}(T)}=O(T^{1/2-\eta_0}),$$
 for some $\eta_0>0$ and the proof is done.
 
 It is definitely possible to optimize these bounds and get an explicit exponent $\eta_0$ in terms of all the constants involved, but this is likely to be non-optimal.

\subsection{Geodesic stretch and bounds on $\delta_0$.} 
In this subsection we derive some lower and upper bounds on $\delta_0$ related to the geodesic stretch of the Teichm\"uller representation $\rho$.
Given $\rho:\Gamma\rightarrow \mathrm{PSL}_2(\R)$ a Teichm\"uller representation, we set
$$I_-(\rho):=\inf_{\gamma \in \mathcal{P}} \frac{\ell_\rho(\gamma)}{\ell(\gamma)},\  I_+(\rho):=\sup_{\gamma \in \mathcal{P}} \frac{\ell_\rho(\gamma)}{\ell(\gamma)}.$$
We prove the following fact.
\begin{propo}
\label{stretch}
Using the above notations, we have
$$\frac{\delta}{2}+\frac{\beta}{2} I_-(\rho) \leq \delta_0\leq \frac{\delta}{2}+\frac{\beta}{2} I_+(\rho).$$
In particular if $I_-(\rho)=I_+(\rho)=1$ then $\delta_0=\frac{\delta+\beta}{2}=1/2+\beta$ which is consistent with the examples computed in $\S 3$.
\end{propo}
\noindent {\it Proof}. We prove only the upper bound, the lower bound being similar. By the variational principle for the topological pressure we have 
$$P(2\beta\psi-2\sigma \tau)=\sup_{\mu} \left ( h_\mu(T)+2\beta \int \psi d\mu-2\sigma \int \tau d\mu\right),$$
where the supremum is taken over all $T$-invariant probability measures. Since this supremum is realized by a unique ergodic (and mixing) equilibrium state, we can replace this supremum
by a sup over all {\it ergodic} $T$-invariant measures. On the other hand, we have 
$$P(\beta\psi-\delta \tau)=\sup_{\mu} \left ( h_\mu(T)+\beta \int \psi d\mu-\delta \int \tau d\mu\right)=0, $$
and we get that
$$ P(2\beta\psi-2\sigma \tau)\leq \sup_{\mu} \left ( \beta \int \psi d\mu+(\delta -2\sigma)\int \tau d\mu\right).$$
Taking 
$$\sigma=\frac{\delta}{2} +\frac{\beta}{2}\sup_\mu \left( \frac{\int \psi d\mu}{\int \tau d\mu} \right),$$
we have $P(2\beta\psi-2\sigma \tau)\leq 0$ and therefore $\delta_0\leq \sigma$. To complete the proof of the upper bound, we just need to show that
$$\sup_\mu \left( \frac{\int \psi d\mu}{\int \tau d\mu} \right)=I_+(\rho).$$
First observe that if $T^nx=x$ is a primitive periodic point of $T$ with $\tau^{(n)}(x)=\ell(\gamma)$ and $\psi^{(n)}(x)=\ell_\rho(\gamma)$, then by taking
$\mu=\frac{1}{n} \sum_{j=0}^{n-1}D_{T^jx}$, where $D_y$ stands for the Dirac measure at $y$, we have
$$\frac{\int \psi d\mu}{\int \tau d\mu}=\frac{\ell_\rho(\gamma)}{\ell(\gamma)}.$$
Since we know that for all $\gamma \in \mathcal{P}$, there exists an $n$-periodic point of $T$ such that both $\tau^{(n)}(x)=\ell(\gamma)$ and $\psi^{(n)}(x)=\ell_\rho(\gamma)$ hold,
we deduce that
$$\sup_\mu \left( \frac{\int \psi d\mu}{\int \tau d\mu} \right)\geq I_+(\rho).$$
On the other hand, if $\mu$ is a $T$ ergodic measure, one can show using Birkhoff's ergodic theorem (see \cite{NaudHP}, Lemma 3.4) that there exists a sequence $(x_p)_{p\in \N}$ of $n_p$-periodic points such that
$$\frac{\int \psi d\mu}{\int \tau d\mu}=\lim_{p\rightarrow \infty}   \frac{\frac{1}{n_p}\psi^{(n_p)}(x_p)}{ \frac{1}{n_p}\tau^{(n_p)}(x_p)}.$$
But we do have for all $p$
$$\frac{\frac{1}{n_p}\psi^{(n_p)}(x_p)}{ \frac{1}{n_p}\tau^{(n_p)}(x_p)}\leq I_+(\rho),$$
and hence the equality is proved. Using \cite{Burger}, we know that if $I_-(\rho)=I_+(\rho)$ then they are both equal to $1$ and the Manhattan curve is a straight line passing through the points $(0,1)$ and $(1,0)$, which implies that $\delta=1+\beta$ and $\delta_0=1/2+\beta$. The proof is done. $\square$

\bigskip

\begin{rem}
The quantities $I_-(\rho)$ and $I_+(\rho)$ are called ``geodesic stretch'' and are related to Thurston's asymetric distance: in particular, we have \cite{PT}
$$\log\left ( \frac{I_+(\rho)}{I_-(\rho)}\right) \leq d_{\mathcal{T}}(X,\rho(\Gamma)\backslash \H),$$
where $d_{\mathcal{T}}(X_1,X_2)$ denotes the usual Teichm\"uller distance related to quasi-conformal deformations. In particular, if the Teichm\"uller distance between
$X$ and the ``deformed'' surface $\rho(\Gamma)\backslash \H$ is small, then $\delta_0$ is close to $(\delta+\beta)/2$.
\end{rem}

\end{section}
\begin{section}{Final remarks}
They are many ways in which this paper can be generalized. It raises also a lot of open questions. We list them below.
\begin{itemize}
\item In this paper we have essentially adressed the case of representations $$\rr:\Gamma\rightarrow \mathrm{SL}_d(\R),$$ 
where $\rr$ is obtained as a composition of Teichm\"uller and the irreducible representation. It would be interesting to consider a more general deformation of that case, which is called the ``Hintchin component'' in the litterature. It is likely that one has to combine techniques from Labourie \cite{Labourie} and the transfer operator methods used in this paper to reach this goal. One could also look at the quasi-Fuchsian deformation of the surface group $\Gamma$ into $\mathrm{PSL}_2(\C)$. We believe
this should follow from our technique without major modifications.

\item By looking at examples, it seems that on the parabola $P_\delta$, $\delta(1-\delta)$ is the only possible eigenvalue. Can one prove that fact in a general setting?
We have shown that if the character of $\rr$ is positive, then $\delta(1-\delta)$ is indeed the bottom of the spectrum of $\Delta_\rr$. Is it a simple eigenvalue? Can one express the base eigenfunction via Patterson-Sullivan theory in higher rank by analogy with the case of geometrically finite hyperbolic manifolds?
\item Based on the examples, we believe that the spectrum of $\Delta_\rr$ always has an essential spectral gap i.e., there exists $\epsilon>0$ such that $\mathrm{Sp}(\Delta_\rr)\cap \mathcal{C}_\delta\setminus \mathcal{C}_{\delta-\epsilon}$ is a finite set. Spectral gaps in the context of negative curvature are often derived from Dolgopyat's techniques \cite{Dolg} and their generalizations to non-compact situations with unitary twists \cite{OW1,MN1,SW1}. Alternatively, when the "trapped set" is a smooth symplectic submanifold and the transverse dynamics are normally hyperbolic, Nonnenmacher and Zworski's micro-local approach \cite{NZ1} provide an explicit essential spectral gap.
In this context of non-unitary representations, Dolgopyat's arguments cannot be applied directly and the micro-local approach is also likely to fail due to the lack of contact/symplectic structure. 
\item The computable examples also suggest that the spectrum of $\Delta_\rr$ may have in some situations a ``band structure'', with the number of bands given by $\mathrm{dim}(V)$. In \cite{FT1}, for Pollicott-Ruelle resonances related to $U(1)$-extensions of Anosov maps, a band structure of the spectrum is rigourously established by Faure and Tsujii, together with some eigenvalue distribution results based on Sj\"ostrand's method \cite{Sjostrand1} for damped wave equations. We also mention the relevant work
of Dyatlov \cite{Dyatlov2} for the existence of band structures in quantum resonances in a setting where the trapped set is normally hyperbolic. It would be interesting to investigate to what extend the above mentioned works could be adapted to our problem.

\item By analogy with the theory of resonances for hyperbolic manifolds, it would be interesting to know if micro-local techniques used by Dyatlov in \cite{Dyatlov1} could be adapted 
to improve our Weyl bounds on the peripheral spectrum of $\Delta_\rr$. The non-smooth nature of Teichm\"uller deformations (Mostow's map is only H\"older) is likely to be an obstacle.
\item There is a current work \cite{PF} on meromorphic continuation of Selberg zeta functions with non-unitary twists, for geometrically finite groups.  Despite the fact that there is no available interpretation yet of zeros as spectra of a twisted Laplacian in this non-compact context, it would be interesting to know if our result can be extended to this framework. 
\item Finally, one can consider hyperbolic manifolds of higher dimension, or more general,
locally symmetric spaces $X=\Gamma\backslash Y$ of real rank one, where $Y$ is a real, complex, or quaternionic hyperbolic space, or the hyperbolic 
Cayley plane. Considering a non-unitary representation of the lattice, the Selberg trace formula, for suitable integral operators, induced by the twisted Laplacian, should be available, and hence the meromorphic continuation of the twisted Selberg zeta function could  follow. Consequently, it would be interesting to see if our results can be extended to this case.
\end{itemize}
\end{section}

\appendix 

\section{Universal bounds for Bergman kernels}
The results below are based on standard facts on Bergman kernels for simply connected domains. However, because we couldn't find the exact statement needed in the classical bet, we include a detailed proof here. In particular, the ``off-diagonal'' estimate is not surprising but seems to be new in this setting.
Let $\mathbb{D}$ denote the unit disc. Let $\Omega\subset \C$ be a bounded, simply connected open set. If $\phi:\Omega \rightarrow \mathbb{D}$ is a conformal representation, 
then one can endow $\Omega$ with a hyperbolic distance $\rho_\Omega(x,y)$ by pulling back the hyperbolic distance $\rho_{\mathbb{D}}(z,z')$ on $\mathbb{D}$:
$$\rho_\Omega(x,y):=\rho_{\mathbb{D}}(\phi(x),\phi(y)).$$
Notice that 
$$\rho_{\mathbb{D}}(z,z')=\inf_{\gamma}\int_0^1 \frac{2\vert \gamma'(t)\vert}{1-\vert \gamma(t)\vert^2}dt, $$
where the infimum runs over all $C^1$ paths $\gamma:[0,1]\rightarrow {\mathbb D}$ such that $\gamma(0)=z$ and $\gamma(1)=z'$. The quantity
$$ \lambda_{\mathbb D}(z):=\frac{2}{1-\vert z\vert^2}$$
is called the density of the hyperbolic metric and we have under the conformal mapping $\phi$ the identity
$$\lambda_\Omega(z)=\frac{2\vert \phi'(z)\vert}{1-\vert \phi(z)\vert^2}.$$
For example, if $\Omega=D(a,r)$ is a euclidean disc centered at $a$ and of radius $r$, then we have
$$\lambda_\Omega(z)= \frac{2r^{-1}}{1-\vert \frac{z-a}{r}\vert^2}.$$
The Bergman Kernel $B_\Omega(z,w)$ is the Schwarz kernel of the orthogonal projection $B:L^2(\Omega)\rightarrow H^2(\Omega)$,
where $H^2(\Omega)$ is as usual the Hilbert space of holomorphic $L^2$ functions on $\Omega$. For the unit disc we have the well known formula \cite[Chapter 1]{Krantz}:
$$B_{\mathbb D}(z,w)=\frac{1}{\pi(1-z\overline{w})^2}.$$
The Bergman kernel has conformal invariance i.e., under the above notations, we have
$$B_\Omega(z,w)=\frac{\phi'(z)\overline{\phi'(w)}}{\pi(1-\phi(z)\overline{\phi(w)})^2}.$$
Using the formula for the hyperbolic distance $\rho_{\mathbb D}(u,v)$ in the disc
$$\cosh^2(\rho_{\mathbb D}(u,v)/2)=\frac{\vert 1-u\overline{v}\vert^2}{(1-\vert u\vert^2)(1-\vert v\vert^2)},$$
we obtain the formula:
\begin{equation}
\label{Berghyp}
\vert B_\Omega(z,w)\vert=\frac{1}{4\pi}\frac{\lambda_\Omega(z)\lambda_\Omega(w)}{\cosh^2(\rho_{\Omega}(z,w)/2)},
\end{equation}
which expresses the kernel in terms of hyperbolic distance and density on $\Omega$.
This formula will be enough to prove the following estimate.
\begin{propo}
\label{Bergman}
Assume that $\Omega$ is a bounded simply connected domain. We have the following facts.
\begin{enumerate}
\item Assume that $z,w\in \Omega$ satisfy $\mathrm{dist}(z,\partial \Omega)\geq c\epsilon$, $\mathrm{dist}(w,\partial \Omega)\geq c\epsilon$, for some $c>0$. Then there exists
$M(c)>0$ such that for all $1\geq \epsilon>0$,
$$\vert B_\Omega(z,w)\vert \leq M(c)\epsilon^{-2}.$$
\item Assume in addition that $z,w\in \Omega$ satisfy: 
$$c\epsilon<\mathrm{dist}(z,\partial \Omega)\leq\epsilon,\ c\epsilon \leq \mathrm{dist}(w,\partial \Omega)\leq \epsilon,\ \vert z-w\vert\geq c\epsilon^\alpha$$ for some $c>0$ and $0<\alpha<1$, then there exists $\widetilde{M}(c)>0$ such that for all $1\geq \epsilon>0$,
$$ \vert B_\Omega(z,w)\vert \leq \widetilde{M}(c)\epsilon^{-1-\alpha}.$$
\end{enumerate}
The above constants $M(c),\widetilde{M}(c)$ do not depend on the domain $\Omega$.

\end{propo}
The first bound is ``a diagonal estimate'' and essentially optimal without additional assumptions on $z,w$. The second bound (``off-diagonal estimate'') shows that provided $z,w$ are sufficiently far away from each other, we can gain over the trivial diagonal bound. These bounds are universal: we do not assume anything on the domain $\Omega$ except
that it is bounded and simply connected.

\bigskip
\noindent {\it Proof}. The hyperbolic density $\lambda_\Omega$ satisfies the domain monotonicity: if $D(z,r)\subset \Omega$ then we have
$$\lambda_\Omega(z)\leq \lambda_{D(z,r)}(z)=2r^{-1}.$$
Using Formula (\ref{Berghyp}) we deduce that under the hypotheses of claim $1$ we have for all $1\geq \epsilon>0$,
$$\vert B_\Omega(z,w)\vert \leq \frac{1}{4\pi}\lambda_\Omega(z)\lambda_\Omega(w)\leq \frac{1}{\pi} \mathrm{dist}(z,\partial \Omega)^{-1}\mathrm{dist}(w,\partial \Omega)^{-1}$$
$$\leq \frac{c^{-2}}{\pi}\epsilon^{-2},$$
and the first bound is proved. 

 The second bound is more subtle: we need to bound from below the hyperbolic distance $\rho_\Omega(z,w)$. To this end, we will use an inequality
due to Beardon \cite{Beardon1}. Given $x,y\in \Omega$, where $\Omega$ is again bounded and simply connected, set
$$\alpha_\Omega(x,y):=\sup_{a,b \in \partial \Omega} \log \left (  \frac{\vert z-a\vert \vert w-b\vert}{\vert z-b\vert \vert w-a\vert}    \right ). $$
This pseudo-metric is called the Apollonian metric by Beardon in \cite{Beardon1}, and satisfies the remarkable bound
$$\alpha_\Omega(x,y)\leq 2\rho_\Omega(x,y). $$
Choosing points $a,b\in \partial \Omega$ such that
$$\vert z-a\vert=\mathrm{dist}(z,\partial \Omega),\  \vert w-b\vert=\mathrm{dist}(w,\partial \Omega),$$
we obtain for $0<\epsilon\leq \epsilon_0(c)$,
$$\alpha_\Omega(z,w)\geq \left \vert \log\mathrm{dist}(z,\partial \Omega)+\log\mathrm{dist}(w,\partial \Omega)-\log\vert z-b\vert-\log \vert w-a\vert \right \vert$$
$$\geq 2\vert \log \epsilon\vert- \vert \log\vert z-b\vert+\log \vert w-a\vert \vert.$$
Using the triangle inequality we have for all $0<\epsilon\leq \epsilon_0(c)$ 
$$\vert z-b\vert\geq \vert z-w\vert-\vert w-b\vert \geq c\epsilon^\alpha-\epsilon\geq \frac{c}{2}\epsilon ^\alpha.$$
Similarly, we obtain also for all $0<\epsilon\leq \epsilon_0(c)$ 
$$ \vert w-a\vert\geq \frac{c}{2}\epsilon ^\alpha,$$
so that we get
$$ \vert \log\vert z-b\vert+\log \vert w-a\vert \vert\leq 2\alpha \vert \log \epsilon \vert +\vert \log(c)\vert.$$
Hence for all $0<\epsilon\leq \epsilon_0(c)$  we have 
$$ (1-\alpha)\vert \log \epsilon \vert -\frac{1}{2}\vert \log c\vert \leq \rho_\Omega(z,w).$$
Going back to formula (\ref{Berghyp}) and plugging the above estimate we get the desired upper bound (valid for all $0<\epsilon\leq \epsilon_0(c)$ )
$$\vert B_\Omega(z,w)\vert \leq \frac{c^{-2}\epsilon^{-2}}{\pi}e^{\frac{1}{2}\vert \log c\vert} \epsilon^{1-\alpha}.$$
The proof is done since for $1\geq \epsilon>\epsilon_0(c)$ we can obviously use the first universal bound. $\square$

We point out that no assumptions on the regularity of $\partial \Omega$ or the convexity of $\Omega$ have been made.

\end{document}